\documentclass[11pt,letterpaper]{article}
\usepackage[T1]{fontenc}
\usepackage[utf8]{inputenc}
\usepackage{lmodern}
\usepackage[margin=1in]{geometry}
\usepackage{amsmath,amssymb,amsthm,mathtools}
\usepackage{graphicx,booktabs,tabularx,enumitem,needspace,float}
\usepackage[ruled,section]{algorithm}
\usepackage{microtype}
\usepackage{xurl}
\usepackage[hidelinks]{hyperref}
\numberwithin{equation}{section}
\newtheorem{theorem}{Theorem}[section]
\newtheorem{lemma}[theorem]{Lemma}
\newtheorem{proposition}[theorem]{Proposition}
\newtheorem{assumption}[theorem]{Assumption}
\theoremstyle{remark}

\newenvironment{keywords}{\par\smallskip\noindent\textbf{Keywords.}\ }{\par\smallskip}
\newenvironment{MSCcodes}{\noindent\textbf{Mathematics Subject Classification (2020).}\ }{\par\medskip}

\setlist{itemsep=2pt,topsep=4pt}
\title{The Query Complexity of Unknown Sparse Gradient Residuals}
\author{Yuriy Dorn\\[0.4em]{\normalsize Lomonosov Moscow State University}}
\date{}
\hypersetup{pdfauthor={Yuriy Dorn}, pdftitle={The Query Complexity of Unknown Sparse Gradient Residuals},
 pdfsubject={Nonconvex finite-sum optimization and query complexity},
 pdfkeywords={nonconvex optimization, finite sums, query complexity, variance reduction}}
\begin{document}
\maketitle
\begin{abstract}
We study optimization of an average of $n$ functions when gradient
predictions are exact except for at most $k$ unknown components.
An initial scan may miss these exceptions.

For a class of smooth problems, we prove that the optimal worst-case
expected number of component-gradient evaluations needed to find an
approximate stationary point is
\[
 \Theta\!\left(\min\left\{nk+\mathcal R,\frac nk\mathcal R\right\}\right),
\]
where $\mathcal R$ represents the optimal cost when the exceptional
components are known. Prediction evaluations are free.

This bound shows that discovering even a few exceptions can cost more
than one full scan and can account for most of the gradient evaluations.
\end{abstract}

\section{Introduction}
A component gradient may be expensive to evaluate even when a useful
approximation is available. If that approximation is exact for most
components, repeated evaluations of those components are unnecessary.
The difficulty is to identify the exceptions. A gradient prediction can
agree with the true gradient at every point queried so far and still differ
elsewhere. Agreement at finitely many points therefore does not establish
that the two gradients agree everywhere.

We consider
\begin{equation}\label{eq:main-objective}
 \min_{x\in\mathbb R^d}F(x),\qquad
 F(x)=\frac1n\sum_{i=1}^n f_i(x)
     =h(x)+\frac1n\sum_{i=1}^n u_i(x).
\end{equation}
For component $i$, a known differentiable function $p_i$ predicts the
gradient through $\nabla p_i$. Write $h=n^{-1}\sum_i p_i$,
$u_i=f_i-p_i$, and $r_i=\nabla u_i$. We call $u_i$ the residual function
and $r_i$ its residual gradient. At most $k$ residual gradients are nonzero
as functions; their coordinates need not be sparse. One query returns
$\nabla f_i(x)$ at a chosen index and point. Subtracting the known
$\nabla p_i(x)$ gives $r_i(x)$. Predictor evaluations are not counted.
Every component retains weight $1/n$, including after an exception is found.

We assume that $F$ is smooth, its initial gap is at most $F_0$, and
\[
 \sum_i\|r_i(x)-r_i(y)\|^2\leq k\delta^2\|x-y\|^2,
 \qquad \sum_i\operatorname{osc}(u_i)\leq2nF_0.
\]
Here $\operatorname{osc}(u_i)=\sup u_i-\inf u_i$. This bound limits the
change in the partial objective when a newly discovered residual function
is added. It is a global assumption in addition to gradient or Hessian
similarity. The individual oscillations and the exceptional indices are
unknown. We compare optimization when these indices are supplied with
optimization when the algorithm must find them.

\Needspace{5\baselineskip}
\subsection{Contributions}
\paragraph{The cost of unknown exceptional indices}
For the parameter range in Theorem~\ref{thm:phase-complexity}, we establish
\[
 \Theta\!\left(\min\left\{nk+\mathcal R,\frac nk\mathcal R\right\}\right),
 \qquad \mathcal R=\frac{k^{3/2}\delta F_0}{n\epsilon^2}.
\]
The theorem covers $\epsilon^2\leq\delta F_0\sqrt{k}/(4n)$.
The terms change order at the scales $k^2$ and $nk$ in $\mathcal R$.
When $n$ is much larger than $k$, these distinguish sampling the population,
finding exceptions, and optimizing their sum. If $n=k$, supplying the
indices does not improve the order. The explicit lower-bound constants
are conservative and do not identify numerical crossover points.
If the exceptional indices are
supplied, the order reduces to $\Theta(\mathcal R)$.
Table~\ref{tab:information-complexity} specifies the information supplied
and the queries counted in each comparison.

\paragraph{Algorithms and stopping rules}
We estimate either the full residual gradient sum or its sum over a growing
set of discovered exceptions, called the working set. After recomputing the
working-set sum exactly, we check the full gradient if the gradient of the
partial objective is small. A failed check adds at least one exception; a
successful check returns a point with gradient norm at most $\epsilon$.
There are at most $k$ failed checks. A common potential is nonincreasing
when the working set expands and satisfies a descent estimate during
optimization. We bound the estimator error up to the random stopping time
to obtain the expected query count. Sampling several differences per
iteration permits stepsizes determined by smoothness alone, while retaining
the same order of query cost. For sampling
from all components or from a supplied exceptional set, we also prove that
the algorithm may stop at the first successful exact gradient check.
The adaptive method in Theorem~\ref{thm:adaptive-working-set} uses the
number of discovered components to choose its expected batch size and
omits differences that an exact recomputation would replace. Its bound
depends on the actual support size $k_*$ and the number $D$ of checks
that discover new components:
\[
 \mathbb E Q\leq(n+2k_*)(\mathbb ED+1)
       +\frac{8k_*\sqrt{k}\delta(F_0+\Omega)}{n\epsilon^2}.
\]
Neither $k_*$ nor $D$ is supplied to the method.

\paragraph{Lower bounds and classes with simpler identification}
We combine an optimization lower bound with smooth functions whose
successive informative gradients belong to components with unknown indices.
Independent copies give the accuracy-dependent discovery bound. An
interrupted-search inequality permits arbitrary interleaving of those
copies, including chains of length one. Convex objectives with bounded
sinusoidal residuals cover the remaining accuracy range. The bounds allow
query points and indices to depend on previous
observations, and the output need not have been queried. For real analytic
residuals, querying each component at one independently sampled point
identifies all exceptions almost surely, giving an $O(n+\mathcal R)$ upper
bound. We also construct predictors with $O(n)$ queries under a
finite-dimensional mixture assumption. Our numerical comparisons include
exceptions that first become visible at successive partial stationary
points. Every competitive baseline may use the exceptional indices that
its own queries identify; supplied-support methods are reported separately.

\subsection{Related work}
Working-set methods alternate optimization over a restricted set of
variables and checks of the full optimality conditions
\cite{strongrules,blitz,blitzws,fireworks}. Our working set records component indices whose residual gradients have
been observed to be nonzero. Those gradients may vanish at earlier query
points.
A bound on the sum of their oscillations controls changes of the partial
objective as this set grows.

PAGE, SPIDER, and MARINA alternate exact gradient recomputations with
updates based on sampled gradient differences
\cite{page,spider,marina}. We use this estimator structure and also account
for queries used to find the exceptional components and check the full
gradient. On a known set of $m$ components, the residual estimator has
mean-square smoothness parameter $\delta\sqrt{mk}/n$; applying the usual
variance-reduction balance gives the optimization term in our upper bound.
Our additional analysis concerns the changing partial objective, the
queries used to discover exceptions, and termination after an exact check.
Similarity-based methods exploit small differences between
component gradient increments~\cite{sliding,svrp,svrs}.
Jiang, Rodomanov, and Stich~\cite{costsaga} study functional similarity
while assigning different costs to arbitrary, random, and delegated
component selection. Their bounds concern communication and local
computation. We count exact component-gradient evaluations. We assume that residual
gradients vanish identically outside a small exceptional set and bound
the sum of the residual oscillations. Predictor evaluations are not counted. These assumptions specify both the function
class and the cost measured by our comparison.

The lower bounds use smooth chains and random rotations
\cite{carmonlowerone,spider}, together with searches for informative
component indices. Han, Xie, and Zhang~\cite{hanxiezhang2024} study
component proximal queries with restrictions on sampling and on the span
of allowable query points. Arjevani et al.~\cite{arjevani2023stochastic}
use searches over random seeds to obtain informative gradient observations.
Our bounds allow arbitrary adaptive component-gradient queries and compare
known and unknown exceptional indices within the same function class.
Recent preprints also study adaptive information acquisition in finite sums.
Zhou~\cite{zhou2026local} assigns successive chain edges to random
components in lower bounds for approximate local minima.
Peng, Tang, and Jia~\cite{peng2026dense} distribute hidden information
across components to obtain lower bounds under individual smoothness.
Our use of unknown component indices is related to these constructions.
The distinction is the comparison between supplied and unknown indices
under exact residual sparsity and a bound on residual oscillations;
neither of those papers states that comparison.
\section{Model and complexity statement}
\label{sec:rb-model}
For each differentiable component $f_i:\mathbb R^d\to\mathbb R$, let
$p_i:\mathbb R^d\to\mathbb R$ be a known differentiable function whose
gradient predicts $\nabla f_i$. Define
\[
 h=\frac1n\sum_i p_i,\qquad u_i=f_i-p_i,\qquad r_i=\nabla u_i,
 \qquad F=h+\frac1n\sum_i u_i.
\]
The exceptional set, also called the support of the residual gradients, is
$\mathcal B=\{i:r_i\not\equiv0\}$. Even for $i\in\mathcal B$, $r_i$
may vanish at all points queried so far. A constant residual function
affects neither the gradient nor the objective gap. All vector
norms are Euclidean, and matrix operator norms are induced by them.

\begin{assumption}[Objective and sparse residuals]\label{ass:ct-objective}
The supplied constants $L,F_0,\delta>0$, $\Omega\geq0$ and $1\leq k\leq n$
satisfy $F(x_0)-\inf F\leq F_0$, $\nabla F$ is $L$-Lipschitz, and
$|\mathcal B|\leq k$. Writing $W_i=\sup u_i-\inf u_i$, assume
\begin{equation}\label{eq:total-potential-budget}
 \frac1n\sum_{i\in\mathcal B}W_i\leq\Omega,
\end{equation}
and, for all $x,y\in\mathbb R^d$,
\begin{equation}\label{eq:gated-rms}
 \sum_{i\in\mathcal B}\|r_i(x)-r_i(y)\|^2
 \leq k\delta^2\|x-y\|^2.
\end{equation}
The individual $W_i$ and the set $\mathcal B$ are not supplied.
\end{assumption}
Condition~\eqref{eq:gated-rms} permits an individual residual gradient
to have Lipschitz constant as large as $\delta\sqrt{k}$.
Condition~\eqref{eq:total-potential-budget} bounds changes of the partial
objective. Both conditions apply to all points. They cannot be established
solely by observing that gradient predictions agree with queried gradients
at finitely many points. Section~\ref{sec:predictor-acquisition} gives
assumptions under which finitely many queries determine the predictors.

\paragraph{Queries and the required output}
A query selects an index $i$ and a point $x$ and returns the exact vector
$\nabla f_i(x)$. We do not count predictor evaluations or other computation.
Evaluating a residual gradient difference at two points costs two queries.
The count includes all component evaluations used to construct predictors
or check the output. We denote the cost of constructing predictors, when
needed, by $C_{\rm cal}$. If a component averages data losses, one query
may require evaluating the gradient on every example in that component.

Algorithms may randomize. They may choose indices and points using all
previous responses, with no bound on the norm of a query point, and return
a point that has never been queried. An admissible algorithm terminates
almost surely and satisfies
$\mathbb E\|\nabla F(\widehat x)\|^2\leq\epsilon^2$ on every instance.
For each algorithm, take the supremum of its expected query count over all
instances and finite dimensions. The infimum of this quantity over admissible
algorithms is the minimax query cost. Our algorithms satisfy the stronger
condition $\|\nabla F(\widehat x)\|\leq\epsilon$, verified by an exact
gradient check before returning. When the exceptional set is supplied before
any query, this check requires only gradients from that set.

\paragraph{Reduction to a common reference function}
The transformed components $\widetilde f_i=h+f_i-p_i=h+u_i$ have the same
average $F$. One query to $f_i$ gives a gradient of $\widetilde f_i$ through
$\nabla\widetilde f_i=\nabla f_i-\nabla p_i+\nabla h$. The reverse
calculation also uses one query. The residual gradients, exceptional set,
smoothness bounds, and oscillations are unchanged. We may therefore use
$f_i=h+u_i$ in the analysis even when the original predictors differ.
\Needspace{22\baselineskip}
\begin{theorem}[Query complexity with unknown and supplied exceptional indices]
\label{thm:phase-complexity}
Let $16\leq k\leq n$ be integers, $\delta,F_0>0$, and
\[
 L_0=\frac{\delta\sqrt{k}}n,\qquad
 L\geq L_0,\qquad c_0=\frac14.
\]
Consider the class in Assumption~\ref{ass:ct-objective} with
$\Omega=2F_0$, using the queries and output requirement defined above.
The reference function is supplied. If constructing it instead costs $O(n)$
queries, this cost does not change the order of the bound for unknown
exceptional indices. For
\[
 0<\epsilon^2\leq c_0\frac{\delta F_0\sqrt{k}}n,\qquad
 \mathcal R=\frac{k^{3/2}\delta F_0}{n\epsilon^2},\qquad
 \mathcal P=\frac{\sqrt{k}\delta F_0}{\epsilon^2}=\frac nk\mathcal R,
\]
the minimax expected query complexity with unknown exceptional indices is
\begin{equation}\label{eq:phase-complexity}
 \Theta\!\left(\min\{nk+\mathcal R,\mathcal P\}\right)
 =\Theta\!\left(\mathcal R+\min\{nk,\mathcal P\}\right).
\end{equation}
The constants are independent of $n,k,L,\delta,F_0,\epsilon$.
The algorithm attaining the upper bound terminates almost surely and
returns a point that passes an exact gradient check.
The lower bound applies to every admissible algorithm defined above;
the dimension of its hard instance may depend on the problem parameters.

If the exceptional set is supplied and there is no requirement to query
all components to construct predictors or check the output, the complexity is
$\Theta(\mathcal R)$. A requirement of $n$ additional component queries changes this
to $\Theta(n+\mathcal R)$.
\end{theorem}

Table~\ref{tab:information-complexity} compares the information available
to the algorithm.
When the exceptional indices of smooth residuals are unknown, the bound
contains the additional term
$\min\{nk,\mathcal P\}$. Since $\mathcal P\geq\mathcal R$, the two
expressions in~\eqref{eq:phase-complexity} agree within a factor of two.
The number of iterations depends on $L$. The choices of recomputation
probabilities and expected batch sizes remove this dependence from the
total query bounds. Sections~\ref{sec:total-potential-budget} and
\ref{sec:lower-bounds} prove the upper and lower bounds, respectively.

\begin{table}[H]
\caption{Query bounds under the parameter range of
Theorem~\ref{thm:phase-complexity}. Supplied indices are known before
any query. The smooth rows give matching lower
and upper bounds; the analytic row gives the upper bound of
Proposition~\ref{prop:analytic-support-probe}. All component-gradient evaluations
are counted.}
\label{tab:information-complexity}
\centering
\begin{tabular}{@{}lll@{}}
\toprule
Residual class & Information and requirements & Query bound\\
\midrule
Smooth & Indices unknown & $\Theta(\min\{nk+\mathcal R,\mathcal P\})$\\
Smooth & Indices supplied & $\Theta(\mathcal R)$\\
Smooth & Indices supplied; $n$ required queries & $\Theta(n+\mathcal R)$\\
Real analytic & Indices unknown & $O(n+\mathcal R)$\\
\bottomrule
\end{tabular}
\end{table}

\paragraph{Accuracy range}
For example, let $n=k^2$ and
$\epsilon^2=\delta F_0\sqrt{k}/(nk)$. Then $\mathcal R=k^2$ and
$\mathcal P=nk=k^3$. As $k$ grows, the unknown-index cost is
$\Theta(k^3)$ and the supplied-index cost is $\Theta(k^2)$.
The additional cost therefore exceeds the $n=k^2$ queries needed for
a single scan.

Put $G=\delta F_0\sqrt{k}/n$. The theorem covers $0<\epsilon^2\leq G/4$,
or equivalently $\mathcal R\geq4k$. The two transition scales are
$\epsilon^2\asymp G/k$ and $\epsilon^2\asymp G/n$, corresponding to
$\mathcal R\asymp k^2$ and $\mathcal R\asymp nk$.
The corresponding intervals fit inside the stated range when $16\leq k<n$.
Their separation is informative when $n$ is much larger than $k$;
at $n=k$, the whole bound has order $\mathcal R$.

The proof uses two kinds of examples. Smooth chains establish the
accuracy-dependent discovery cost at sufficiently small $\epsilon$.
Bounded sinusoidal residuals establish linear query costs for the remaining
accuracies. Combining the two gives the stated uniform order, with explicit
constants in Section~\ref{sec:lower-bounds}. Those constants remain
conservative: the theorem does not locate numerical crossover points for
a fixed moderate $n$ and $k$. Since $\mathcal P\geq4n$ and $nk\geq n$,
adding $O(n)$ queries to construct predictors does not change the order.
\section{Upper bounds with bounded residual oscillation}
\label{sec:total-potential-budget}

We work under Assumption~\ref{ass:ct-objective}. If the predictors must be
constructed from queries, their construction costs $C_{\rm cal}$. This cost
is added once to each bound in this section.

\subsection{A working-set method with exact gradient checks}
The working set $S$ contains components whose residual gradients have been observed
to be nonzero. For $S\subseteq\mathcal B$, define the partial objective
and its residual gradient sum by
\[
 \Phi_S=h+\frac1n\sum_{i\in S}u_i,\qquad
 R_S=\frac1n\sum_{i\in S}r_i,\qquad
 \mathcal H=\frac{k\delta}{n},\qquad
 \Gamma=F_0+\Omega.
\]
Each partial objective $\Phi_S$ is $(L+\mathcal H)$-smooth. Its gradient
differs from $\nabla F$ by the sum of the omitted residual gradients,
divided by $n$. By Cauchy--Schwarz and~\eqref{eq:gated-rms}, this difference
has Lipschitz constant at most $\mathcal H$.
Choose
\begin{equation}\label{eq:aggregate-gate-parameters}
 \begin{aligned}
 \eta&=\frac1{2(L+\mathcal H)},&
 p&=\frac{2\eta\mathcal H}{\sqrt k},\\
 \tau&=\eta\mathcal H\sqrt k,&
 b&=\lceil\tau\rceil,\qquad s=\frac\tau b.
 \end{aligned}
\end{equation}
These parameters satisfy $0<p\leq1$ and $0<s\leq1$.
At each visit, we recompute $R_S(x)$ exactly with probability $p$.
After each executed step, we estimate its change. When $S\ne\varnothing$,
this estimate uses $\tau$ differences of residual gradients on average.
The batch size $b$ allows this expected count to exceed one.

\begin{algorithm}[H]
\caption{Working-set method with exact gradient checks}\label{alg:siopt-gate}
\begin{minipage}{\linewidth}
\textbf{Input:} $x_0,\epsilon$ and $(\eta,p,\tau,b,s)$ from
\eqref{eq:aggregate-gate-parameters}; the known gradient $\nabla h$ and
access to queries of the residual gradients.
Initialize $S=\varnothing$, $x=x_0$, and $\bar v=0$.
A visit consists of one pass through the following operations. The last
visit may return a point before taking a step.
\begin{enumerate}
\item Independently of the past, recompute $R_S(x)$ with probability $p$
by querying every $i\in S$, and set $v=R_S(x)$. If no recomputation is
made, set $v=\bar v$.
\item If step 1 recomputed $R_S(x)$ and
$\|\nabla h(x)+v\|\leq\epsilon/2$, query all components.
Return $x$ if $\|\nabla F(x)\|\leq\epsilon$. Otherwise add every
$i\notin S$ with $r_i(x)\ne0$ to $S$ and reset $v=R_S(x)$ using
the responses from this check.
\item Set $x^+=x-\eta(\nabla h(x)+v)$.
If $S=\varnothing$, set $\bar v^+=0$ without a query.
Otherwise, after fixing $x^+$, draw mutually independent
$Z_j\sim\operatorname{Bernoulli}(s)$ and
$I_j\sim\operatorname{Uniform}(S)$, $j=1,\ldots,b$.
For each $j$ with $Z_j=1$, query component $I_j$ at $x$ and $x^+$.
Set
\begin{equation}\label{eq:gated-difference}
 \bar v^+=v+\frac1b\sum_{j=1}^b\frac{Z_j}{s}\frac{|S|}{n}
                     \bigl(r_{I_j}(x^+)-r_{I_j}(x)\bigr).
\end{equation}
If every $Z_j$ is zero, retain $v$ without querying residual gradients. Set
$(x,\bar v)=(x^+,\bar v^+)$.
\end{enumerate}
\end{minipage}
\end{algorithm}
Recomputing the sum over $S$ costs at most $k$ queries. Each sampled
difference costs two queries, and checking the full gradient costs $n$.
When $S\ne\varnothing$, the expected cost of the sampled differences is
$2bs=2\tau$. Indices remain in the working set after they are added.
The tests for a zero residual gradient and for the final gradient norm
use exact component-gradient responses.

\begin{theorem}[Query bound for the working-set method]
\label{prop:gated-total-budget}
Under Assumption~\ref{ass:ct-objective}, Algorithm~\ref{alg:siopt-gate} terminates almost surely and returns a point satisfying
$\|\nabla F(x)\|\leq\epsilon$. If $V$ counts all visits, including the
terminal visit, and $Q$ counts all queries, then
\begin{align}
 \mathbb E V&\leq\frac{8\Gamma}{\eta\epsilon^2}+\frac{k+1}{p},
 \label{eq:aggregate-gate-visits}\\
 \mathbb E Q&\leq(n+2k)(k+1)
                  +\frac{32k^{3/2}\delta\Gamma}{n\epsilon^2}.
 \label{eq:aggregate-gate-cost}
\end{align}
In particular, if $\Omega\leq2F_0$ and
$\mathcal R=k^{3/2}\delta F_0/(n\epsilon^2)$, then
$\mathbb E Q\leq(n+2k)(k+1)+96\mathcal R$.
\end{theorem}

\begin{proof}
We first restrict each sum to a fixed finite number of visits. After the
algorithm terminates, we set all further movements and estimator errors to
zero and include no further descent terms. Whether a visit is reached is
determined before deciding whether to recompute the residual gradient
sum. These finite sums let us bound the expected stopping time without assuming in
advance that it is finite.

\paragraph{A potential for all working sets}
Put $W_i=\operatorname{osc}(u_i)$. Add constants to the residual functions
so that $0\leq u_i\leq W_i$, and subtract their average from $h$.
These additive constants affect no gradient or objective gap. Define
\[
 \Psi_S(x)=\Phi_S(x)
                +\frac1n\sum_{i\in\mathcal B\setminus S}W_i.
\]
Then $\Psi_S\geq F\geq\inf F$ and
$\Psi_\varnothing(x_0)-\inf F\leq\Gamma$.
Adding a set $A$ at a fixed point changes this potential by
$n^{-1}\sum_{i\in A}(u_i(x)-W_i)\leq0$.
Adding components therefore contributes no increase when we sum the
descent inequalities. The initial difference remains bounded by $\Gamma$.

\paragraph{Estimator error up to termination}
Let $\bar e=\bar v-R_S(x)$ before the recomputation decision, and let
$e=v-R_{S^{\rm new}}(x)$ after that decision and any addition to $S$.
An exact recomputation gives $e=0$, including on the last visit.
If the sum is not recomputed, then $e=\bar e$. New components are added
only when their contributions have been evaluated exactly. Conditional on
the history before this decision,
$\mathbb E\|e\|^2=(1-p)\|\bar e\|^2$.
For fixed $x$ and $x^+$, the sampled sum added to $v$
in~\eqref{eq:gated-difference} is an unbiased estimate of the change in $R_S$. Independence of its $b$ terms bounds its conditional variance by
\[
 \frac{|S|}{bsn^2}\sum_{i\in S}\|r_i(x^+)-r_i(x)\|^2
       \leq\frac{\mathcal H^2}{\tau}\|x^+-x\|^2.
\]
The indicators $Z_j$ are independent; using a single indicator for the
whole batch would not give the factor $1/b$. Subtracting the true change in $R_S$ from the sampled sum gives an error
whose inner product with $e$ has conditional expectation zero. Writing
$a_t=\mathbb E\|e_t\|^2$ and taking the movement to be zero when no
step executes gives $a_0=0$ and
\[
 a_{t+1}\leq(1-p)a_t+
       \frac{(1-p)\mathcal H^2}{\tau}
                         \mathbb E\|x^{t+1}-x^t\|^2.
\]
Summing this recursion over the fixed finite number of visits yields
\begin{equation}\label{eq:aggregate-stopped-error}
 \sum_t\mathbb E\|e_t\|^2
 \leq C\sum_t\mathbb E\|x^{t+1}-x^t\|^2,
 \qquad C=\frac{(1-p)\mathcal H^2}{p\tau}
                 \leq\frac1{2\eta^2}.
\end{equation}
The recursion is conditioned only on the history available before each draw.

\paragraph{Descent and the number of visits}
For an executed step, smoothness and the identity
$x^+-x=-\eta(\nabla\Phi_S(x)+e)$ give
\[
 \Phi_S(x^+)\leq\Phi_S(x)
 -\frac\eta2\|\nabla\Phi_S(x)\|^2
 -\left(\frac1{2\eta}-\frac{L+\mathcal H}{2}\right)\|x^+-x\|^2
 +\frac\eta2\|e\|^2.
\]
Sum these inequalities over the executed steps, using $\Psi_S$ to account
for additions to $S$, and apply~\eqref{eq:aggregate-stopped-error}. The
coefficient multiplying the sum of squared movements is nonnegative because
$1/(2\eta)-(L+\mathcal H)/2-\eta C/2\geq0$.
Consequently
\begin{equation}\label{eq:aggregate-gradient-energy}
 \mathbb E\sum_{\text{executed visits}}
       \|\nabla\Phi_{S^{\rm new}}(x)\|^2\leq\frac{2\Gamma}{\eta}.
\end{equation}
If the full gradient check fails, at least one new exceptional component
is added. Otherwise all omitted residual gradients would vanish at that point, and
the gradient of the partial objective would equal the full gradient. There are
therefore at most $k$ failed checks.

Separate visits according to the exact gradient norm of the partial
objective before the recomputation decision. If
$\|\nabla\Phi_S(x)\|>\epsilon/2$, no full gradient check can occur,
$S$ stays unchanged, and a step is taken. Equation~\eqref{eq:aggregate-gradient-energy}
bounds the expected number of these visits by $8\Gamma/(\eta\epsilon^2)$.
Let $V_{\rm low}$ count visits with $\|\nabla\Phi_S(x)\|\leq\epsilon/2$.
Each such visit leads to a full gradient check with probability $p$.
There are at most $k+1$ full checks on every run, including the last one.
Both arrival at a visit and the gradient threshold condition are determined
before the independent recomputation decision. Hence
$p\mathbb E V_{\rm low}\leq k+1$ for every finite limit on the visit count.
Adding the two bounds and applying monotone convergence
proves~\eqref{eq:aggregate-gate-visits} and almost-sure termination.
The algorithm returns only after checking the exact gradient norm.

\paragraph{Query count}
Full gradient checks cost at most $n(k+1)$. Conditional on reaching a
visit, exact recomputation uses at most $kp$ queries in expectation. Each
executed step uses at most $2bs=2\tau$ queries in expectation for differences
of residual gradients, because the samples are drawn after $x$ and $x^+$ are fixed.
Therefore
\[
 \mathbb E Q\leq n(k+1)+(kp+2\tau)\mathbb E V
 \leq n(k+1)+4\eta\mathcal H\sqrt k
           \left(\frac{8\Gamma}{\eta\epsilon^2}+\frac{k+1}{p}\right),
\]
which is~\eqref{eq:aggregate-gate-cost}.
\end{proof}

\paragraph{Effect of batching on the visit bound}
The stepsize in~\eqref{eq:aggregate-gate-parameters} depends on the
smoothness of the partial objective. To use at most one sampled difference
per step, we may choose
$\eta_{\rm one}=\min\{1/[2(L+\mathcal H)],1/(2\mathcal H\sqrt k)\}$.
This gives the same query bound. Allowing several differences improves the
corresponding bound on visits by
the factor $\eta/\eta_{\rm one}
=\max\{1,\mathcal H\sqrt k/(L+\mathcal H)\}\leq\sqrt k$.
This compares the proved bounds on visits. The numbers observed in a run
also depend on the objective and random samples.

\subsection{Checking omitted components one at a time}
\label{sec:sequential-complement}
Suppose an exact recomputation at $x$ gives
$g_S=\nabla\Phi_S(x)$ with $\|g_S\|\leq\epsilon/2$.
The check of all components can be replaced as follows. If $|S|=k$, return $x$.
Otherwise, query $[n]\setminus S$ at $x$ in a fresh uniformly random
order. At the first nonzero residual gradient, add its index to $S$, reset
$v=R_S(x)$ using the recomputed sum and the new response, and resume
optimization.
If all components outside $S$ have zero residual gradient at $x$, return $x$.
Every query is counted. A zero response at $x$ does not permanently
exclude the component from later checks.

Returning immediately when $|S|=k$ is valid: $S\subseteq\mathcal B$ and
$|\mathcal B|\leq k$ then imply $S=\mathcal B$. This stopping condition cannot occur if the exceptional set has fewer than
$k$ members.

\begin{proposition}[Query bound for sequential gradient checks]
\label{prop:sequential-complement}
Replace the full gradient checks in Theorem~\ref{prop:gated-total-budget}
by this procedure. The method terminates almost surely with
$\|\nabla F(x)\|\leq\epsilon/2$ and retains the visit
bound~\eqref{eq:aggregate-gate-visits}. Its expected query count satisfies
\begin{equation}\label{eq:sequential-total-cost}
 \mathbb E Q\leq C_{\mathrm{scan}}+2k(k+1)
          +\frac{32k^{3/2}\delta\Gamma}{n\epsilon^2},
\end{equation}
where $C_{\mathrm{scan}}\leq n(k+1)$ is the expected total number of
queries to components outside $S$. A scan queries these components in
random order until the first nonzero response, or until all have been queried.
Conditional on the history when a scan starts, put
\[
 N=|[n]\setminus S|,\qquad
 a=|\{i\notin S:r_i(x)\ne0\}|.
\]
The expected query count of that scan is
\begin{equation}\label{eq:first-hit-cost}
 \begin{cases}
 (N+1)/(a+1),&a\geq1,\\
 N,&a=0.
 \end{cases}
\end{equation}
Returning when $|S|=k$ requires no queries outside $S$.
\end{proposition}

\begin{proof}
Each interrupted scan adds one new member of $\mathcal B$.
If the scan finishes without a nonzero response, then
$g_S=\nabla F(x)$. The same equality holds when $|S|=k$. At most $k$
scans add a component, and at most one further check returns the final point.

Each addition to $S$ follows an exact recomputation and introduces no
estimator error. It cannot increase the potential $\Psi_S$ from
Theorem~\ref{prop:gated-total-budget}. That proof therefore gives the same
error recursion, bound on the sum of squared gradient norms, and count of
visits with gradient norm above $\epsilon/2$. At a visit with gradient
norm at most $\epsilon/2$, a check occurs with conditional probability $p$;
there are at most $k+1$ checks. This proves the same visit bound and
almost-sure termination. Exact recomputations and sampled differences cost
at most $(kp+2\tau)\mathbb E V$. Substituting the visit bound and adding
the queries outside $S$ proves~\eqref{eq:sequential-total-cost}.

For $a\geq1$, the positions with nonzero responses divide the $N-a$
zero-response positions into $a+1$ gaps with the same expected length.
The expected first gap is $(N-a)/(a+1)$; the first nonzero response adds
one query. If $a=0$, all $N$ queries
are required. These conditional expectations may be summed over all
scans that occur, even though earlier responses determine when a new
scan starts.
\end{proof}

\paragraph{The cost of stopping a scan early}
A complete scan can discover several components at once. Stopping at the
first nonzero response discovers one component and may require further
optimization before the next scan.
Consequently~\eqref{eq:first-hit-cost} does not imply a smaller total cost.
For example, suppose exactly $b$ exceptional components have nonzero
residual gradients at every point checked until all are discovered,
starting with $S=\varnothing$. The $b$ interrupted
scans have expected total cost
\[
 \sum_{j=0}^{b-1}\frac{n-j+1}{b-j+1}
 =b+(n-b)(H_{b+1}-1),\qquad H_q=\sum_{j=1}^q\frac1j.
\]
A final scan confirming that all omitted residual gradients are zero
adds $n-b$ queries when $b<k$. It is unnecessary when $b=k$. A complete
initial scan identifies all $b$ indices for $n$ queries. Conversely,
when exactly one omitted residual gradient is nonzero at the point
being checked, stopping at the first nonzero response discovers the
same index for $(N+1)/2$ expected queries. These counts cover only the
queries used to discover components. The total cost also includes the
optimization steps between checks.

\subsection{Adapting the work to the discovered components}
\label{sec:adaptive-working-set}
The preceding methods use a batch size determined by the supplied upper
bound $k$, even when few components have been discovered. They also sample
a difference before learning whether the next exact recomputation will
replace it. We can avoid both costs. The resulting bound depends on the
number of checks that discover new components and on the actual support
size, which may be smaller than $k$.

Let $\overline L$ be a supplied smoothness bound for every partial
objective the method can form. The choice $\overline L=L+\mathcal H$
always suffices under Assumption~\ref{ass:ct-objective}; a smaller bound
may follow from additional structure. Choose
\[
 0<\eta\leq\frac1{2\overline L},\qquad
 p=\min\left\{1,\frac{2\eta\delta\sqrt{k}}n\right\}.
\]
For a nonempty working set of size $m$, use
\begin{equation}\label{eq:adaptive-working-batch}
 \tau_m=\frac{\eta\delta\sqrt{k}}n m,
 \qquad b_m=\lceil\tau_m\rceil,\qquad s_m=\frac{\tau_m}{b_m}.
\end{equation}
When $p<1$, $\tau_m=pm/2$. When $p=1$, every sum is recomputed and
no differences are sampled. No differences are sampled when $m=0$.
Algorithm~\ref{alg:adaptive-working-set} draws the recomputation decision
before paying for a difference. It uses the same tolerance $\epsilon$
for the exact partial gradient and the full gradient: if the latter check
fails, at least one omitted residual gradient must be nonzero.

\begin{algorithm}[H]
\caption{Adaptive working-set method}\label{alg:adaptive-working-set}
\begin{minipage}{\linewidth}
\textbf{Input:} $x_0,\epsilon$, the preceding parameters $\eta,p$, the
known gradient $\nabla h$, and residual-gradient queries.
Initialize $S=\varnothing$, $x=x_0$, $x_{\rm prev}=x_0$, and $v_{\rm prev}=0$.
At each visit:
\begin{enumerate}
\item Put $m=|S|$ and draw an independent Bernoulli variable with success
probability $p$. On success, query every component in $S$ at $x$ and set
$v=R_S(x)$. On failure, if $m=0$ set $v=0$; otherwise draw independent
$Z_j\sim\operatorname{Bernoulli}(s_m)$ and
$I_j\sim\operatorname{Uniform}(S)$ for $j=1,\ldots,b_m$, and set
\[
 v=v_{\rm prev}+\frac1{b_m}\sum_{j=1}^{b_m}
       \frac{Z_j}{s_m}\frac mn
       \bigl(r_{I_j}(x)-r_{I_j}(x_{\rm prev})\bigr).
\]
Query both endpoints only for terms with $Z_j=1$.
\item If step 1 was an exact recomputation and
$\|\nabla h(x)+v\|\leq\epsilon$, do the following.
If $m=k$, return $x$.
Otherwise query every component outside $S$ and put
$w=n^{-1}\sum_{i\notin S}r_i(x)$.
Return $x$ if $\|\nabla h(x)+v+w\|\leq\epsilon$.
If the check fails, add all omitted indices with $r_i(x)\ne0$ to $S$
and replace $v$ by $v+w$.
\item Store $x_{\rm prev}=x$ and $v_{\rm prev}=v$, then set
$x=x-\eta(\nabla h(x)+v)$ and begin the next visit.
\end{enumerate}
\end{minipage}
\end{algorithm}

\Needspace{21\baselineskip}
\begin{theorem}[Cost in terms of discoveries and support size]
\label{thm:adaptive-working-set}
Under Assumption~\ref{ass:ct-objective} and the preceding choices of
$\overline L,\eta,p$, Algorithm~\ref{alg:adaptive-working-set}
terminates almost surely and returns a point with exact gradient norm at
most $\epsilon$. Put $k_*=|\mathcal B|$. Let $D$ be the number of failed
full gradient checks, and let $m_j$ be the working-set size immediately
before check $j$, for $j=0,\ldots,D$, including the terminal check. Then
$D\leq k_*$, and the total number of queries to omitted components is
\[
 C=\sum_{j=0}^{D}(n-m_j){\bf1}_{\{m_j<k\}}.
\]
The numbers $V$ of visits and $Q$ of queries satisfy
\begin{align}
 \mathbb E V&\leq\frac{2\Gamma}{\eta\epsilon^2}
                         +\frac{\mathbb E D+1}{p},
 \label{eq:adaptive-working-visits}\\
 \mathbb E Q&\leq
 \mathbb E\left[C+(2-p)\sum_{j=0}^{D}m_j\right]
 +\frac{4(2-p)k_*\sqrt{k}\delta\Gamma}{n\epsilon^2}.
 \label{eq:adaptive-working-cost}
\end{align}
In particular,
\begin{equation}\label{eq:adaptive-working-simple}
 \mathbb E Q\leq(n+2k_*)(\mathbb E D+1)
                  +\frac{8k_*\sqrt{k}\delta\Gamma}{n\epsilon^2}.
\end{equation}
No value of $D$, $k_*$, or the individual residual oscillations is needed
by the method.
\end{theorem}

The first term in~\eqref{eq:adaptive-working-cost} charges the components
actually queried at checks. In particular, a check can discover many
components at once. The second term scales with $k_*$ rather than the
supplied support bound $k$. Setting $D\leq k_*\leq k$ and
$\Omega\leq2F_0$ gives the worst-case bound
$(n+2k)(k+1)+24\mathcal R$. The order of the minimax bound is unchanged.

For $k_*>0$, write the residual smoothness bound as
$k_*\delta_*^2=k\delta^2$, where $\delta_*=\delta\sqrt{k/k_*}$.
The optimization term with supplied support then scales as
$k_*^{3/2}\delta_*=k_*\sqrt{k}\delta$.
The adaptive method obtains this dependence without knowing $k_*$,
while also paying for discovery and replacing $F_0$ by $F_0+\Omega$.
This statement concerns the working-set bound; selection between it and
population sampling still uses the supplied parameters.

\begin{proof}
The set before a visit is exactly the set used for the preceding step.
All components added at the preceding visit were included in the exact
estimate before that step. Thus the deferred difference estimates the
change of one fixed residual sum, even when the preceding visit enlarged
the set. Both endpoints and this set are fixed before the current
recomputation decision and difference samples.

For $m>0$, the conditional variance of the sampled difference is at most
\[
 \frac{m}{n^2\tau_m}\sum_{i\in S}
       \|r_i(x^+)-r_i(x)\|^2
 \leq\frac{\delta\sqrt{k}}{\eta n}\|x^+-x\|^2.
\]
For $m=0$, the error is zero. With $a_t$ denoting the expected squared
estimator error after any additions at visit $t$, independence and
unbiasedness give
\[
 a_{t+1}\leq(1-p)a_t+
       \frac{(1-p)\delta\sqrt{k}}{\eta n}
                        \mathbb E\|x^{t+1}-x^t\|^2,
 \qquad a_0=0.
\]
As before, first sum over a fixed finite number of visits and set errors
and movements after termination to zero. If $p<1$, the sum of expected
squared errors is at most $(1-p)/(2\eta^2)$ times the sum of expected
squared movements, because $p=2\eta\delta\sqrt{k}/n$. If $p=1$, every
error is zero. Use $\overline L$ in the smoothness inequality and
$\eta\leq1/(2\overline L)$ in its movement coefficient. The common potential
$\Psi_S$ therefore gives~\eqref{eq:aggregate-gradient-energy} without
change.

Classify a reached visit before drawing its recomputation decision.
If the exact partial gradient norm exceeds $\epsilon$, no check or set
enlargement occurs. The gradient-energy bound limits the expected number
of these visits to $2\Gamma/(\eta\epsilon^2)$. At any other visit, a
check occurs with probability $p$. A failed check adds a new component,
so there are at most $k_*+1$ checks. Finite truncation and monotone
convergence first prove almost-sure termination. Summing the conditional
check probabilities then gives $p\mathbb E V_{\rm low}=\mathbb E D+1$,
which proves~\eqref{eq:adaptive-working-visits}. Returning when $m=k$
is valid because $S\subseteq\mathcal B$ and $|\mathcal B|\leq k$.
Every other return uses the full gradient obtained from the complementary
queries and the exact sum on $S$.

Conditional on reaching a visit with working-set size $m$, updating the
residual estimate costs at most
\[
pm+2(1-p)\tau_m=(2-p)pm
\]
queries in expectation. The identity holds for $p<1$ by
$\tau_m=pm/2$, and for $p=1$ because no differences are sampled.
This also holds at the initial visit, where no
difference is sampled. On visits with partial gradient norm above
$\epsilon$, use $m\leq k_*$, the preceding bound on their number, and
$p/\eta\leq2\delta\sqrt{k}/n$.
On the remaining visits, both their classification and $m$ are determined
before the recomputation decision. Consequently
\[
 p\mathbb E\sum_{\text{low visits}}|S|
       =\mathbb E\sum_{j=0}^{D}m_j.
\]
Adding the $C$ complementary queries proves
\eqref{eq:adaptive-working-cost}. Finally, $C\leq n(D+1)$ and
$\sum_jm_j\leq k_*(D+1)$ give~\eqref{eq:adaptive-working-simple}.
\end{proof}

\paragraph{Comparison with exact working-set gradients}
An exact-gradient version uses the same discovery rule and
$x^+=x-\nabla\Phi_S(x)/\overline L$.
The usual smoothness inequality gives a decrease of at least
$\|\nabla\Phi_S(x)\|^2/(2\overline L)$ at each executed step.
After a failed full check, the enlarged partial gradient equals the full
gradient and has norm greater than $\epsilon$. Thus every executed step
has partial gradient norm greater than $\epsilon$.
Using the common potential, the step count is at most
$2\overline L\Gamma/\epsilon^2$. Its query count is therefore at most
\[
 C+k_*+\frac{2k_*\overline L\Gamma}{\epsilon^2},
\]
where $C$ is its own complementary-scan count. The two optimization terms
compare the residual scale $\delta\sqrt{k}/n$ with the partial-objective
smoothness $\overline L$. Small batches alone do not imply a lower total
cost when the exact-gradient method can take larger steps.

\subsection{Estimation over all components or a known support}
The same estimator can sample all components or a known set containing
every exceptional component. Let $A\subseteq[n]$
be a known set containing $\mathcal B$, put $m=|A|\geq2$, and define
\[
 R=\frac1n\sum_{i\in A}r_i,\qquad
 \ell_A=\frac{\delta\sqrt{mk}}n.
\]
Choose
\[
 0<\eta_A\leq\frac1{2L},\qquad
 p_A=\min\left\{1,\frac{2\eta_A\ell_A}{\sqrt m}\right\},\qquad
 \tau_A=\eta_A\ell_A\sqrt m,
\]
and set $b_A=\lceil\tau_A\rceil$ and $s_A=\tau_A/b_A$.
The largest allowed stepsize is $\eta_A=1/(2L)$. The smaller choice
$\eta_A=\min\{1/(2L),1/(\ell_A\sqrt m)\}$ has $\tau_A\leq1$ and
$b_A=1$, so each step samples at most one difference. If $p_A=1$, the sum
is recomputed at every step and no differences are sampled.

Query every component in $A$ at $x_0$ to initialize $v_0=R(x_0)$, and return $x_0$ if
$\|\nabla h(x_0)+v_0\|\leq\epsilon$.
Otherwise execute $x_{t+1}=x_t-\eta_A(\nabla h(x_t)+v_t)$.
After computing $x_{t+1}$, independently recompute $R(x_{t+1})$ with
probability $p_A$. In this case, query every component in $A$, set
$v_{t+1}=R(x_{t+1})$, and return $x_{t+1}$
if $\|\nabla h(x_{t+1})+v_{t+1}\|\leq\epsilon$.
If no recomputation is made, draw mutually independent
$Z_j\sim\operatorname{Bernoulli}(s_A)$ and
$I_j\sim\operatorname{Uniform}(A)$, $j=1,\ldots,b_A$, and set
\[
 v_{t+1}=v_t+\frac1{b_A}\sum_{j=1}^{b_A}\frac{Z_j}{s_A}\frac mn
                   \bigl(r_{I_j}(x_{t+1})-r_{I_j}(x_t)\bigr).
\]
For each $j$ with $Z_j=1$, query component $I_j$ at $x_t$ and $x_{t+1}$.
If every $Z_j$ is zero, retain $v_t$. Continue until an exact recomputation
shows that the gradient norm is at most $\epsilon$.
The choice $A=[n]$ requires no support information. If the exceptional
set is supplied, take it as $A$ and, if necessary, add known
indices whose residual gradients vanish identically to make $m\geq2$; for $k\geq2$ this keeps $m\leq k$.

\begin{proposition}[Query bounds with and without supplied support]
\label{prop:tracking-population}
Suppose $2\leq k\leq n$ and Assumption~\ref{ass:ct-objective} holds,
except that the oscillation condition~\eqref{eq:total-potential-budget}
is not required. The preceding method terminates almost surely and returns
a point whose exact gradient norm has been checked to be at most $\epsilon$. Its number $T_A$ of executed steps satisfies
\[
 \mathbb E T_A\leq\frac{2F_0}{\eta_A\epsilon^2}+\frac1{p_A}.
\]
Its expected query cost, including initialization, satisfies
\begin{equation}\label{eq:tracking-population-cost}
 \mathbb E Q\leq3m+\frac{8m\sqrt k\delta F_0}{n\epsilon^2}.
\end{equation}
Sampling all components costs $O(n+\mathcal P)$, where
$\mathcal P=\sqrt k\delta F_0/\epsilon^2$. Sampling a supplied support
costs $O(k+\mathcal R)$. The final gradient checks query $n$ and at most
$k$ components, respectively. If the protocol also requires queries to all
$n$ components, their cost is added to the bound with supplied support.
\end{proposition}

\begin{proof}
Because $A$ contains every component with a nonzero residual gradient, $\nabla h+R=\nabla F$.
If the initial gradient check succeeds, the claims hold with $T_A=0$ and $Q=m$.
Otherwise $\|\nabla F(x_0)\|>\epsilon$.
First sum over at most $N$ steps, where $N$ is fixed, and set movements
and estimator errors after termination to zero. Independence of the
sampled terms and their inclusion indicators bounds the conditional
variance of the difference estimate by
\[
 \frac{m}{b_As_A n^2}\sum_{i\in A}\|r_i(x_{t+1})-r_i(x_t)\|^2
       \leq\frac{\ell_A^2}{\tau_A}\|x_{t+1}-x_t\|^2.
\]
For $e_t=v_t-R(x_t)$, every exact recomputation makes the error zero,
including the recomputation at termination. If the sum is not recomputed,
the algorithm continues. The error in the sampled difference has inner
product with $e_t$ of conditional expectation zero.
Consequently $a_t=\mathbb E\|e_t\|^2$ satisfies $a_0=0$ and
\[
 a_{t+1}\leq(1-p_A)a_t+
 \frac{(1-p_A)\ell_A^2}{\tau_A}
                   \mathbb E\|x_{t+1}-x_t\|^2.
\]
Summing this recursion over the first $N$ steps gives
\[
 \sum_{t<N}\mathbb E\|e_t\|^2
 \leq C_A\sum_{t<N}\mathbb E\|x_{t+1}-x_t\|^2,
 \qquad
 C_A=\frac{(1-p_A)\ell_A^2}{p_A\tau_A}\leq\frac1{2\eta_A^2}.
\]
For $p_A<1$, use $p_A\tau_A=2\eta_A^2\ell_A^2$; for $p_A=1$,
the errors and $C_A$ vanish. At each executed step, smoothness yields
\[
 F(x_{t+1})\leq F(x_t)
 -\frac{\eta_A}2\|\nabla F(x_t)\|^2
 -\left(\frac1{2\eta_A}-\frac L2\right)\|x_{t+1}-x_t\|^2
 +\frac{\eta_A}2\|e_t\|^2.
\]
Let $J_N=\min\{T_A,N\}$. Sum the descent inequalities through $x_{J_N}$,
use $F\geq\inf F$, and apply the error bound. Since
$\eta_A\leq1/(2L)$, the coefficient multiplying the sum of squared
movements is nonnegative. Hence
\[
 \mathbb E\sum_{t<J_N}\|\nabla F(x_t)\|^2\leq\frac{2F_0}{\eta_A}.
\]

Consider the point reached after each executed step, before deciding
whether to recompute the residual gradient sum. Let $L_N$ count such
points with true gradient norm at most $\epsilon$, and let $H_N$ count those with
larger norm. At a point counted in $L_N$, the algorithm terminates with
conditional probability $p_A$. It can terminate only once, so
$p_A\mathbb E L_N\leq1$. The number of steps that end with gradient
norm above $\epsilon$ minus the number that start above this threshold is
\[
 H_N-\sum_{t<J_N}{\bf1}_{\{\|\nabla F(x_t)\|>\epsilon\}}
 ={\bf1}_{\{\|\nabla F(x_{J_N})\|>\epsilon\}}-1\leq0,
\]
because $\|\nabla F(x_0)\|>\epsilon$. The bound on the sum of squared
gradient norms therefore gives
$\mathbb E H_N\leq2F_0/(\eta_A\epsilon^2)$.
Since $J_N=L_N+H_N$, monotone convergence proves the step bound and
almost-sure termination. Every returned point passes an exact gradient check.

Initialization costs $m$ queries. Conditional on reaching a step, the
expected cost of updating the residual estimate is
$c_A=mp_A+2(1-p_A)\tau_A$. Whether that step is reached is determined
before the recomputation decision and samples are drawn. Therefore
\[
 \mathbb E Q\leq m+c_A\mathbb E T_A.
\]
For $p_A<1$, $\tau_A/p_A=m/2$ and
$c_A=2\eta_A\ell_A\sqrt m(2-p_A)$. Substitution gives
\[
 \mathbb E Q\leq m(3-p_A)+
       \frac{4(2-p_A)\ell_A\sqrt m F_0}{\epsilon^2}
 \leq3m+\frac{8\ell_A\sqrt m F_0}{\epsilon^2}.
\]
For $p_A=1$, $c_A=m$ and $m/\eta_A\leq2\ell_A\sqrt m$, giving
the stronger bound $2m+4\ell_A\sqrt m F_0/\epsilon^2$.
Using $\ell_A\sqrt m=m\sqrt k\delta/n$ proves
\eqref{eq:tracking-population-cost}.
\end{proof}

\paragraph{Using support information acquired during population sampling}
If queries identify $k$ distinct components with nonzero residual gradients,
they identify the complete exceptional set. At the next exact recomputation,
the population method may restrict its sampling set to those indices.
The recomputation resets the estimator exactly. Its probability
$p_A=\min\{1,2\eta_A\delta\sqrt{k}/n\}$ does not depend on the sampling
set size. The conditional variance bound and the proof of the stopping-time
bound are unchanged, while the cost per step can only decrease.
Thus the population query bound remains valid with this change.

\paragraph{Selecting an upper bound}
When $\Omega\leq2F_0$, select the smaller of the valid bounds
\[
 U_{\rm gate}=(n+2k)(k+1)+24\mathcal R,
 \qquad
 U_{\rm pop}=3n+8\mathcal P.
\]
Here $U_{\rm gate}$ uses Theorem~\ref{thm:adaptive-working-set}.
Under Assumption~\ref{ass:ct-objective}, both methods check the exact
gradient norm before returning.
This yields $O(n+\min\{nk+\mathcal R,\mathcal P\})$ queries.
A predictor-construction cost of $O(n)$ leaves this order unchanged when
$\mathcal P\gtrsim n$. The smoothness bound $L$ affects the number of
optimization steps, but does not appear in either query term that depends
on $\epsilon$.
\section{Lower bounds}
\label{sec:lower-bounds}
The proof uses chain constructions and bounded sinusoidal residuals.
When the exceptional indices
are known, independent chains give $\Omega(\mathcal R)$. When the component
assigned to each link is unknown, a single chain gives $\Omega(nk)$ at fine
accuracy. Several such chains in orthogonal coordinate blocks give
$\Omega(\mathcal P)$ at coarser accuracy: each unfinished chain contributes
to the squared gradient norm. The sinusoidal examples cover the remaining
coarse-accuracy range with linear lower bounds in $n$ and $k$.

\subsection{A bounded smooth chain}
\label{sec:siopt-chain}
The accuracy-dependent constructions use the smooth chain of Carmon et al.
\cite{carmonlowerone}. A chain is a sum of links, each depending on at
most two successive coordinates. We first define the chain and prove
the gradient bound used below. Let
\[
 \Psi(t)=\begin{cases}
 0,&t\leq1/2,\\
 \exp\bigl(1-(2t-1)^{-2}\bigr),&t>1/2,
 \end{cases}
 \qquad
 \Phi(t)=\sqrt e\int_{-\infty}^{t}e^{-s^2/2}\,ds.
\]
For $z\in\mathbb R^K$ put
\[
 \begin{aligned}
 d_1(z)&=-\Phi(z_1),\qquad c_K(z)=\sum_{j=1}^K d_j(z),\\
 d_j(z)&=\Psi(-z_{j-1})\Phi(-z_j)-\Psi(z_{j-1})\Phi(z_j),\quad j\geq2.
 \end{aligned}
\]
The function $\Psi$ and all its derivatives vanish when $t\leq1/2$.
Thus a link with preceding coordinate of magnitude below $1/2$ has zero
value and gradient in a neighborhood. A query at such a point reveals
no information about that link. The radial map below sends all query
points into a bounded ball. The added quadratic term keeps the gradient
norm large outside the part of the domain used by the construction.

\begin{lemma}[Bounds for the chain and its gradient]
\label{lem:siopt-chain}
For $K\geq1$, the links satisfy $|d_j|<12$,
$\operatorname{osc}(d_j)\leq24$, $\|\nabla d_j\|\leq23$ and
$\|\nabla^2d_j\|_{\rm op}\leq410$. Moreover,
\[
 c_K(0)-\inf c_K\leq12K,\qquad
 \|\nabla c_K\|\leq23\sqrt K,\qquad
 \|\nabla^2c_K\|_{\rm op}<154.
\]
Let $U$ have $K$ orthonormal columns, set $R=230\sqrt K$, and define
\[
 \rho(y)=\frac{y}{\sqrt{1+\|y\|^2/R^2}},\qquad
 q_U(y)=c_K(U^\top\rho(y)),\qquad H_U(y)=q_U(y)+\|y\|^2/10.
\]
Then $\|D\rho\|_{\rm op}\leq1$,
$\|D^2\rho[v,w]\|\leq6\|v\|\|w\|/R$, and both $q_U,H_U$ have
Hessian norm at most $512$. If the last coordinate of $U^\top\rho(y)$
has magnitude below $1/2$, then $\|\nabla H_U(y)\|>1/2$.
\end{lemma}
\begin{proof}
The scalar factors satisfy $0\leq\Psi<e$, $\Psi(1)=1$,
$0<\Phi<\sqrt{2\pi e}$, $|\Phi'|\leq\sqrt e$, $|\Phi''|\leq1$,
and $|\Psi'|\leq\sqrt{54/e}$. At most one of $\Psi(t),\Psi(-t)$
is nonzero. On $t>1/2$, writing $a=(2t-1)^{-2}$ gives
$\Psi''(t)=8e e^{-a}a^2(2a-3)$; its nonzero extrema at
$a=(9\pm\sqrt{33})/4$ have magnitude below $33$.
Each gradient coordinate of the chain is consequently bounded by
$\sqrt{108\pi}+e^{3/2}<23$, and each Hessian row sum by
$33\sqrt{2\pi e}+2\sqrt{54}+e<154$.
For an individual link, the squared gradient norm is at most
$108\pi+e^3<23^2$. We use the looser Hessian bound $410$ when combining
the bounds for several links.
The value bound follows from $e\sqrt{2\pi e}<12$.
Since $c_K(0)=-\Phi(0)<0$ and $c_K\geq-12K$, the gap bound follows.
Direct differentiation proves the two bounds for $\rho$. The chain rule
then bounds the Hessians of $q_U,H_U$ by $154+0.6$ and $154+0.8$.

To prove the lower bound on the gradient norm, put $z=U^\top\rho(y)$, $r=\|y\|$, and
$a=(1+r^2/R^2)^{-1/2}$. Every coordinate of $g=\nabla c_K(z)$ is
nonpositive. If $|z_K|<1/2$, let $j$ be the first index with $|z_j|<1$.
The link ending at this coordinate has derivative at most $-1$: for $j=1$ this uses
$\Phi'(z_1)>1$, and otherwise $|z_{j-1}|\geq1$, $\Psi(1)=1$,
and $\Phi'(\pm z_j)>1$. The next link can only decrease this derivative.
Writing $v_j$ for column $j$ of $U$ gives
\[
 \nabla H_U(y)=aUg-\frac{a^3}{R^2}y\langle y,Ug\rangle+\frac y5.
\]
If $r\leq R/2$, then $a\geq2/\sqrt5$ and
$|\langle v_j,y\rangle|=|z_j|/a<1/a\leq\sqrt5/2$. The radial correction
in that direction is at most $a^2r\|g\|/R^2\leq1/20$.
The component of the sum of the three terms in this direction has magnitude at
least $2/\sqrt5-1/20-\sqrt5/10>1/2$.
If $r\geq R/2$, use the triangle inequality and
$\|D\rho\|_{\rm op}\leq a$ to obtain
$\|\nabla H_U(y)\|\geq r/5-23a\sqrt K
\geq(23-46/\sqrt5)\sqrt K>1/2$.
\end{proof}
\subsection{An accuracy lower bound when the exceptional indices are known}
\label{sec:accuracy-lower}
We first bound the queries needed to achieve the required accuracy when
the exceptional indices are supplied. We place independent smooth chains
in known coordinate blocks. Their directions are unknown, so the algorithm
must query each block repeatedly. The construction follows the finite-sum
chain argument in~\cite{spider}. We prove directly that it applies to
adaptively chosen query points and indices.

\begin{theorem}[Accuracy lower bound for bounded residual gradients]
\label{thm:bounded-residual-lower}
Let $2\leq k\leq n$, $\delta,F_0,\epsilon>0$, $C=512$, and
\[
 A:=\frac{\delta F_0\sqrt{k}}{768Cn\epsilon^2}\geq4.
\]
There is a finite-dimensional family with a known quadratic reference $h$
proportional to the squared Euclidean norm, predictors $p_i=h$,
and a known exceptional
set $\mathcal B$ of size $k$ such that $F(0)-\inf F\leq F_0$, the objective is
$L_0$-smooth with $L_0=\delta\sqrt{k}/n$, nonexceptional residual gradients vanish, and
\[
 \frac1k\sum_{i\in\mathcal B}\|r_i(x)-r_i(y)\|^2
 \leq\delta^2\|x-y\|^2,
 \qquad
 \sup_{i,x}\|r_i(x)\|^2
 \leq\frac{529}{12C}\frac{\delta nF_0}{\sqrt{k}}.
\]
For every randomized algorithm using exact component-gradient queries
that terminates almost surely and guarantees $\mathbb E\|\nabla F(\widehat x)\|^2\leq\epsilon^2$
on every member of the family, some member has an expected query count satisfying
\begin{equation}\label{eq:bounded-residual-lower-cost}
 \mathbb E N\geq2^{-23}
 \frac{k^{3/2}\delta F_0}{n\epsilon^2}.
\end{equation}
The exceptional residual functions also satisfy
$\operatorname{osc}(u_i)\leq2nF_0/k$, so this family belongs to the
objective class of Theorem~\ref{thm:phase-complexity}.
\end{theorem}
\label{app:bounded-residual-lower}

\begin{proof}[Proof of Theorem~\ref{thm:bounded-residual-lower}]
Use the chain and radial map of Section~\ref{sec:siopt-chain}. Set
\[
 K=2\lfloor A/2\rfloor,\qquad R=230\sqrt K,\qquad
 \alpha=\frac{nF_0}{12kK},\qquad
 \ell_0^2=\frac{C\alpha}{\delta\sqrt k},\qquad
 Q=\lfloor kK/4\rfloor.
\]
Choose the integer dimension of each block to satisfy
\begin{equation}\label{eq:accuracy-block-dimension}
 m\geq K+1+8R^2\log\bigl(16kK(Q+1)\bigr),\qquad d=km.
\end{equation}
Let $C_i:\mathbb R^m\to\mathbb R^d$, $i\in[k]$, be known isometric
coordinate embeddings with mutually orthogonal ranges. Independently draw
$U_i=[v_{i,1},\ldots,v_{i,K}]$ uniformly from the $m$-by-$K$
matrices with orthonormal columns. We call such matrices frames. Define
\[
 \begin{aligned}
 h(x)&=\frac{\alpha}{10n\ell_0^2}\|x\|^2,& f_i&=h+u_i,\\
 u_i(x)&=\alpha q_{U_i}(C_i^\top x/\ell_0)\quad(i\leq k),&
 u_i&=0\quad(i>k).
 \end{aligned}
\]
Thus $\mathcal B=[k]$ is known and
$F(x)=n^{-1}\sum_{i=1}^k\alpha H_{U_i}(C_i^\top x/\ell_0)$.
All quantities in this construction other than the random frames are known
to the algorithm.

\paragraph{Verifying the objective assumptions}
By orthogonality and Lemma~\ref{lem:siopt-chain},
\[
 \frac1k\sum_{i=1}^k\|r_i(x)-r_i(y)\|^2
 \leq\frac{C^2\alpha^2}{k\ell_0^4}\|x-y\|^2
 =\delta^2\|x-y\|^2.
\]
The block-diagonal objective Hessian has norm at most
$C\alpha/(n\ell_0^2)=\delta\sqrt k/n$. The initial objective gap is at most
$12k\alpha K/n=F_0$, and
\[
 \|r_i(x)\|^2\leq529K\alpha^2/\ell_0^2
 =\frac{529}{12C}\frac{\delta nF_0}{\sqrt k}.
\]
Each link has absolute value below $12$. Composing the chain with the radial
map cannot enlarge its range. Hence
$\operatorname{osc}(u_i)\leq24K\alpha=2nF_0/k$ globally.

\paragraph{Comparing actual and auxiliary responses}
Fix the algorithm's random seed. For the proof, define an auxiliary oracle
that has revealed the first $s_i$ columns in each block, initially none.
On a query to $i\leq k$, it reveals the next column if $s_i<K$ and
increases $s_i$ by one.
It then returns the gradient of $h$ plus that of the truncated residual function
\[
 \alpha\sum_{j=1}^{s_i}
 d_j\bigl(U_i^\top\rho(C_i^\top x/\ell_0)\bigr).
\]
This expression uses only the revealed columns. The radius in $\rho$
remains $R=230\sqrt K$, regardless of how many columns have been revealed.
A query to $i>k$ receives the gradient of the known reference. The algorithm
uses only the returned gradients; the additional column information is
used in the analysis. Define the auxiliary process on all outcomes. We next
construct an event on which every returned gradient equals the actual one.

Condition on the auxiliary history, including all revealed columns,
\emph{before} a query. Every unrevealed column is uniform in the orthogonal
complement of the columns already revealed in its block. This space has
dimension at least $m-K$. Applying the radial map to the query in that
block gives a vector of norm at most $R$, so
\[
 \Pr\!\left\{
 |\langle v_{i,j},\rho(C_i^\top x/\ell_0)\rangle|\geq\tfrac12
 \mid\text{preceding auxiliary history}\right\}
 \leq2\exp\!\left(-\frac{m-K-1}{8R^2}\right).
\]
The same estimate holds for every unrevealed column at an arbitrary
output measurable from the final history of the auxiliary process. There are at most
$kK(Q+1)$ tests over $Q$ queries and that output.
Equation~\eqref{eq:accuracy-block-dimension} therefore bounds the probability
of any failed test by $1/8$. This bound does not condition the distribution
of the frame on earlier successful tests.

If all tests succeed, every omitted link vanishes in a neighborhood of
each query. The first link omitted after a column is revealed depends on
the newly revealed coordinate. That coordinate was still unknown and was
tested \emph{before} the query. The preceding coordinates of all later
links were also tested. Thus the truncated and actual gradients agree.
Induction shows that the actual and auxiliary processes have identical
queries, responses, and outputs on this event. This argument
permits arbitrary adaptive indices, unbounded query points, and an unqueried
output. Each query reveals at most one column in one block. After at most
$Q$ calls, at least $k/2$ blocks still have unrevealed columns. Their last
columns are therefore unrevealed, and the corresponding coordinates have
magnitude below one half at the output. Lemma~\ref{lem:siopt-chain} and
orthogonality then give, with probability at least $7/8$,
\[
 \|\nabla F(\widehat x)\|^2
 \geq\frac{k\alpha^2}{8n^2\ell_0^2}
 =\frac{\delta F_0\sqrt k}{96CnK}
 \geq8\epsilon^2.
\]
Averaging over the algorithm's seed gives the same conclusion for any
algorithm limited to this number of queries.

\paragraph{Algorithms with a random stopping time}
Run a general stopping algorithm until it returns or is about to request
its $(Q+1)$st gradient; in the latter case return zero. Denote this truncated
output by $\widehat x_Q$. It agrees with the original output $\widehat x$
on $\{N\leq Q\}$.
Averaging over the family drawn before the queries and the algorithm's
seed, the preceding bound and the assumed bound on the expected squared
gradient norm imply
\[
 \begin{aligned}
 \frac78&\leq\Pr\{\|\nabla F(\widehat x_Q)\|^2\geq8\epsilon^2\}\\
 &\leq\Pr(N>Q)+\Pr\{\|\nabla F(\widehat x)\|^2\geq8\epsilon^2\}\\
 &\leq\Pr(N>Q)+\frac18.
 \end{aligned}
\]
Therefore $\mathbb E N\geq3Q/4$. Since $Q\geq kK/8$ and $K\geq A/2$,
\[
 \mathbb E N\geq\frac{3kA}{64}
 =\frac{k^{3/2}\delta F_0}{16384Cn\epsilon^2}.
\]
With $C=512$, this proves~\eqref{eq:bounded-residual-lower-cost} for some
member of the family. The query count includes any queries used to
construct the predictors.
\end{proof}
\subsection{The cost of discovering exceptional components}
\label{sec:gated-discovery-lower}
We assign successive links of a smooth chain to components whose indices
are unknown. The next direction is revealed only when the algorithm queries
the component holding that link. We account both for adaptively chosen
query points and for the search over component indices.

\begin{theorem}[Lower bound for discovering the exceptional components]
\label{thm:gated-discovery-lower}
Let $16\leq k\leq n$, $\delta,F_0>0$, and $C_*=1024$. Define
\[
 B^2=\frac{\delta nF_0}{\sqrt{k}},\qquad
 W=\frac{2nF_0}{k},\qquad L_0=\frac{\delta\sqrt{k}}n,
 \qquad
 \epsilon_0^2=\frac{\delta F_0}{768C_*n\sqrt{k}}.
\]
There is a finite-dimensional family with a known quadratic reference
proportional to the squared Euclidean norm, predictors $p_i=h$,
and exactly $k$ unknown
exceptional components. Every member satisfies $F(0)-\inf F\leq F_0$ and is
$L_0$-smooth. Its residual functions $u_i$ and residual gradients $r_i$ satisfy
\[
 \sup_{i,x}\|r_i(x)\|\leq B,\qquad
 \operatorname{osc}(u_i)\leq W,\qquad
 \frac1k\sum_{i\in\mathcal B}\|r_i(x)-r_i(y)\|^2
      \leq\delta^2\|x-y\|^2.
\]
For every $0<\epsilon\leq\epsilon_0$ and every randomized algorithm using
exact component-gradient queries that terminates almost surely and guarantees
$\mathbb E\|\nabla F(\widehat x)\|^2\leq\epsilon^2$ on every family member,
some member has an expected query count satisfying
\begin{equation}\label{eq:gated-discovery-lower}
 \mathbb E N\geq\frac{13}{2048}nk.
\end{equation}
The algorithm may choose arbitrary component indices and query points, and may
return an unqueried point. Queries used to construct predictors are included.
The exceptional indices are not supplied.
\end{theorem}
\label{app:gated-discovery-lower}

\begin{proof}[Proof of Theorem~\ref{thm:gated-discovery-lower}]
We draw the objective before any queries, verify its assumptions, and bound
the queries needed to find the component assigned to each successive link.
An auxiliary process lets us compare this search with the actual oracle.

\paragraph{Assigning chain links to unknown component indices}
Use the links $d_1,\ldots,d_k$ and radial map $\rho$ of
Section~\ref{sec:siopt-chain} with chain length $K=k$. Set
\[
 R=230\sqrt{k},\qquad \alpha=\frac{nF_0}{12k},\qquad
 \ell^2=\frac{C_*\alpha}{\delta\sqrt{k}}.
\]
Draw a uniformly random orthonormal frame
$U=[v_1,\ldots,v_k]\in\mathbb R^{d\times k}$ and, independently, a uniformly
random ordered list $(\lambda_1,\ldots,\lambda_k)$ of distinct indices in $[n]$.
Define
\begin{align*}
 h(x)&=\frac{\alpha}{10n\ell^2}\|x\|^2,&
 u_{\lambda_j}(x)&=\alpha d_j(U^\top\rho(x/\ell)),\\
 u_i(x)&=0\quad(i\notin\{\lambda_1,\ldots,\lambda_k\}),&
 f_i&=h+u_i.
\end{align*}
The complete objective is exactly
\[
 F(x)=\frac\alpha n
       \left(c_k(U^\top\rho(x/\ell))+\frac1{10}\|x/\ell\|^2\right).
\]
All quantities other than the random frame and the ordered assignment of
component indices are known to the algorithm.

\paragraph{Bounds on residual gradients and the objective}
Lemma~\ref{lem:siopt-chain} gives
$|d_j|<12$, $\operatorname{osc}(d_j)\leq24$, $\|\nabla d_j\|\leq23$
and $\|\nabla^2d_j\|_{\rm op}\leq410$.
Each link depends on at most two adjacent
coordinates. Hence for every vector $v$,
\[
 \sum_j\|\nabla^2d_j(z)v\|^2\leq2\cdot410^2\|v\|^2.
\]
Let $q_j(y)=d_j(U^\top\rho(y))$. Since $\|D\rho\|_{\rm op}\leq1$ and
$\|D^2\rho[v,w]\|\leq6\|v\|\|w\|/R$, the chain rule and Minkowski's
inequality for the Euclidean norm of the concatenated Hessian-vector
products give
\[
 \left(\sum_j\|\nabla^2q_j(y)v\|^2\right)^{1/2}
 \leq\left(\sqrt2\,410+\frac{6\cdot23\sqrt{k}}R\right)\|v\|
 <C_*\|v\|.
\]
Integrating along a line segment and substituting the chosen scale gives
the stated root-mean-square smoothness bound for the residual gradients. Also
\[
 \|r_i(x)\|^2\leq\frac{529\alpha^2}{\ell^2}
 =\frac{529}{12C_*}\frac{\delta nF_0}{\sqrt{k}}\leq B^2,
 \qquad \operatorname{osc}(u_i)\leq24\alpha=W.
\]
By Lemma~\ref{lem:siopt-chain}, the chain composed with the radial map,
plus the quadratic term, has initial gap at most $12k$ and Hessian norm
at most $512$.
Thus $F(0)-\inf F\leq F_0$ and its smoothness is at most
$512\alpha/(n\ell^2)=\delta\sqrt{k}/(2n)\leq L_0$.
All these bounds hold globally. At the origin only the first link has a nonzero residual gradient.

\paragraph{Revealing successive links with an auxiliary oracle}
Fix the algorithm and its random seed. Let $m$ be the number of successive
links whose component indices and frame columns have been revealed,
initially zero. At query $(i,x)$, if $i=\lambda_j$ for some $j\leq m$,
return the exact gradient. If $m<k$ and $i=\lambda_{m+1}$, reveal
that index and $v_{m+1}$, increase $m$ by one, and return the exact
gradient. This revelation occurs regardless of whether the gradient
itself determines the new direction. For every other component index
not yet revealed, return only $\nabla h(x)$.

Define this process on all outcomes. Its replies depend only on the
columns already revealed. Conditional on the queries, responses, and revealed columns, the
unrevealed columns therefore form a uniformly distributed orthonormal
frame in the orthogonal complement of the revealed columns. The assignment of component indices is independent
of the frame and preserves this conditional distribution.
The search argument allows the algorithm to use all revealed indices
and columns. To compare with the actual oracle, run the given algorithm
while withholding these additional revelations. The comparison below bounds
the probability that the two processes return different gradients.

\paragraph{Agreement at the queries and at the output}
At any auxiliary query, condition on the preceding history. Each currently
unrevealed frame column is a uniform unit vector in a subspace of dimension
at least $d-k$, and $\|\rho(x/\ell)\|\leq R$. The concentration bound for a uniformly distributed unit vector gives
\[
 \Pr\!\left\{|\langle v_j,\rho(x/\ell)\rangle|\geq\tfrac12
             \mid\hbox{preceding auxiliary history}\right\}
 \leq2\exp\!\left(-\frac{d-k-1}{8R^2}\right).
\]
The same statement applies to an arbitrary final output measurable from
that history, before any further direction is revealed. A sufficient dimension
is
\[
 d\geq k+1+8R^2\log(32k(Q+1)),\qquad Q=\lfloor nk/64\rfloor.
\]
A union bound over at most $k(Q+1)$ tests, covering all queries and the
output, bounds the failure probability by $1/16$. We evaluate this
probability in the auxiliary process without conditioning the frame
distribution on the projection bounds at earlier queries.

If all tests succeed, every auxiliary response equals the actual gradient.
It remains only to check an index $\lambda_j$, $j\geq m+2$, assigned to
a later link.
Its preceding coordinate has magnitude less than $1/2$, so $\Psi$ and all its
derivatives vanish there. Its link is identically zero near the query point,
and its true residual gradient is zero. Queries to revealed indices and to the next index return exact
gradients by definition. Consequently an ordinary
algorithm and its auxiliary simulation have identical queries, responses,
and outputs whenever all tests succeed. The frame and component assignment
were drawn before any query, so the objective is fixed throughout.

\paragraph{Cost of finding each successive component index}
At the start of stage $j$, the next index is uniform among $n-j+1$ remaining
indices. A failed query within that stage excludes only the queried index from
being the next one. Previous stages impose no further ordering restrictions
on the remaining indices once their successful indices are known. A query to a
component already revealed, or excluded during the current stage, cannot
reveal the next link.

Fix the history at the start of a stage and the algorithm's randomness.
Consider the order in which it would query distinct remaining indices if
every query to an unrevealed component returned the gradient of the known
reference. Responses at revealed components remain unchanged. Append any indices
it would not query, in a fixed order. This order is independent of the next
unknown index. Until that index is queried, these hypothetical responses
and the actual auxiliary responses coincide. The rank of the next index is conditionally uniform
on $\{1,\ldots,n-j+1\}$. After the index is found, the remaining
assignment is again uniform. By induction, the total query count needed
to reveal every link stochastically dominates
\[
 D=\sum_{j=1}^k U_j,\qquad
 U_j\ \hbox{independent},\quad
 U_j\sim\mathrm{Uniform}\{1,\ldots,n-j+1\}.
\]
This remains true with arbitrary adaptive point choices, repeated indices,
and the algorithm's randomization. Each stage still requires finding its
component index even when all previously revealed directions are available.

For $j\leq\lceil k/2\rceil$, the range size $n-j+1$ is at least $n/2$.
Comparing the exponential sum to its integral gives
\[
 \mathbb E e^{-2U_j/n}
 \leq\frac{1-e^{-2(n-j+1)/n}}{2(n-j+1)/n}\leq1-e^{-1}.
\]
Thus
\[
 \Pr(D\leq Q)
 \leq e^{2Q/n}(1-e^{-1})^{\lceil k/2\rceil}
 \leq e^{k/32}(1-e^{-1})^{k/2}<\frac1{16},\qquad k\geq16.
\]
With probability at least $7/8$, both processes have the same output and the
last frame column remains unrevealed after at most $Q$ queries.

\paragraph{The gradient bound and random stopping times}
The gradient bound in Lemma~\ref{lem:siopt-chain} applies whenever the
last coordinate of the chain argument has magnitude below $1/2$. Thus the
common output of the two processes on the event above satisfies
\[
 \|\nabla F(\widehat x)\|^2>
 \frac{\alpha^2}{4n^2\ell^2}
 =\frac{\delta F_0}{48C_*n\sqrt{k}}
 \geq16\epsilon^2.
\]
The output test was already included in the union bound, so the output need
not have been queried.

For an arbitrary almost-surely terminating algorithm, write $\widehat x$
for its output. Run it until return or immediately before a requested
$(Q+1)$st gradient, returning zero in the latter case; denote this output
by $\widehat x_Q$. They agree on $\{N\leq Q\}$.
Averaging over the objective family and the algorithm's seed, the bound
for at most this many queries and the assumed bound on the expected
squared gradient norm give
\[
 \begin{aligned}
 \frac78&\leq\Pr\{\|\nabla F(\widehat x_Q)\|^2\geq16\epsilon^2\}\\
 &\leq\Pr(N>Q)+\Pr\{\|\nabla F(\widehat x)\|^2\geq16\epsilon^2\}\\
 &\leq\Pr(N>Q)+\frac1{16}.
 \end{aligned}
\]
Therefore $\mathbb E N\geq13Q/16\geq13nk/2048$, since $nk\geq256$.
Averaging implies that some fixed family member attains this expected lower
bound. This count includes every query used to construct the predictors.
\end{proof}

\subsection{Search cost under adaptive allocation}
\label{sec:adaptive-search-cost}
The following estimate permits chains of any positive length. It accounts
for an algorithm that switches between searches and stops unsuccessful
searches early.

\begin{lemma}[Cost of an interrupted search]\label{lem:interrupted-search}
Suppose a stage has one successful index, uniform among $M$ candidates
conditional on the history when the stage begins. Before success, responses
from this stage give no information except whether the queried index is
successful. Randomness used in other searches is independent of that index
conditional on this history. For any adaptive allocation of a finite query
budget, let $Z$ count distinct candidates queried during this stage and let
$J$ indicate that the stage succeeds. Then
\[
 \mathbb E Z\geq\frac{M+1}{2}\mathbb E J.
\]
For $m$ chains, independently draw a uniform ordered list of $q\leq N$
distinct successful indices from $N$ candidates in each chain. At each
stage use the same response rule: before its successful index is queried,
responses give no information about the remaining order beyond excluding
the current query from being that stage's successful index.
Let $C_Q$ be the number of completed chains
after at most $Q$ queries. Arbitrary adaptive interleaving is allowed. Then
\begin{equation}\label{eq:completed-chain-budget}
 \mathbb E C_Q\leq\frac{4Q}{q(2N-q+3)}\leq\frac{4Q}{Nq}.
\end{equation}
\end{lemma}
\begin{proof}
Condition on entry into a stage and expose the independent randomness of
the algorithm and the other searches. Consider its continuation if every
candidate tested in this stage fails. Let $t\leq M$ be the number of
distinct candidates it would test before the global budget expires, and
append untested candidates in a fixed order. The successful index has a
uniform rank $U$ in this order, independent of this continuation. The
actual and hypothetical histories coincide until that index is tested.
Thus $J=\mathbf 1_{\{U\leq t\}}$ and $Z=\min\{U,t\}$. Consequently
\[
 \mathbb E[Z\mid t]
 =\frac{t(2M-t+1)}{2M}
 \geq\frac{M+1}{2}\frac{t}{M}
 =\frac{M+1}{2}\mathbb E[J\mid t].
\]
Repeated queries only increase the paid cost. If a stage is never entered,
both variables are zero, so the inequality also holds without conditioning
on entry. For distinct successive labels, failures in previous stages impose
no restrictions on the remaining order after their successful labels have
been revealed. Hence stage $j$ has $M=N-j+1$ candidates.

Apply the first inequality to every stage in each chain. Every completed
chain has completed all its stages. The queried candidates in these
stages use disjoint portions of the total budget. Therefore
\[
 \begin{aligned}
 Q&\geq\sum_{a=1}^m\sum_{j=1}^{q}\frac{N-j+2}{2}
       \Pr\{\text{stage }(a,j)\text{ completes}\}\\
  &\geq\frac{q(2N-q+3)}4\mathbb E C_Q
   \geq\frac{Nq}{4}\mathbb E C_Q.
 \end{aligned}
\]
The conditional argument includes information obtained in other chains;
it does not condition on the set of chains that eventually complete.
\end{proof}

\subsection{Lower bounds at coarse accuracy}
\label{sec:coarse-accuracy-lower}
A one-dimensional family gives an explicit lower bound when the target
accuracy is a fixed fraction of the objective's smoothness times its
initial gap. The residuals in this family are analytic and the objective
is convex. Thus these bounds do not require the delayed revelation used
by the longer smooth chains.

\begin{proposition}[Coarse-accuracy lower bounds]\label{prop:coarse-lower}
Let $1\leq k\leq n$, $G=\delta F_0\sqrt{k}/n$, and
$0<\epsilon^2<G/2$. Under Assumption~\ref{ass:ct-objective} with $\Omega=2F_0$ and
$L\geq\delta\sqrt{k}/n$, the worst-case expected query cost of
an admissible algorithm is at least
\begin{equation}\label{eq:coarse-lower}
 \begin{aligned}
 &\frac{n+1}{2}\left(1-\frac{2\epsilon^2}{G}\right)
 &&\text{for unknown exceptional indices},\\
 &k\left(1-\frac{2\epsilon^2}{G}\right)
 &&\text{for supplied exceptional indices}.
 \end{aligned}
\end{equation}
The bounds hold for convex objectives with analytic residual functions.
\end{proposition}
\begin{proof}
Put $L_0=\delta\sqrt{k}/n$. First choose an index $\lambda$ uniformly
in $[n]$ and, independently, a sign $\sigma$ uniformly in $\{-1,1\}$.
In one dimension define
\[
 \begin{gathered}
 A=nF_0,\qquad b^2=\frac{2A}{\delta\sqrt{k}},\qquad
 h(x)=\frac{L_0x^2}{4},\qquad
 \\
 u_\lambda(x)=-\sigma A\sin(x/b),\qquad u_i=0\ (i\ne\lambda).
 \end{gathered}
\]
The objective Hessian is $(L_0/2)(1+\sigma\sin(x/b))$, which lies in
$[0,L_0]$. The initial objective gap is at most $A/n=F_0$, the total
residual oscillation is $2A=2nF_0$, and the residual gradient has
Lipschitz constant $A/b^2=\delta\sqrt{k}/2$. Thus every member belongs
to the stated class.

We may give the algorithm the sign as additional information whenever it
queries the exceptional component, including at a point where the queried
residual gradient vanishes. Before this event all responses equal
$\nabla h$, so the sign remains independent of its observations and its
output. Write $a=A/(nb)$ and $y=x/b$. Since $L_0b/2=a$,
\[
 \mathbb E_\sigma|F'(x)|^2
 =a^2\bigl(y^2+\cos^2y\bigr)\geq a^2=\frac G2,
\]
where $\sin^2y\leq y^2$. This inequality holds at every point. It
therefore applies conditionally at any output returned before the
exceptional index is queried. The required expected squared gradient norm
implies that the probability of querying this index is at least
$1-2\epsilon^2/G$. Apply Lemma~\ref{lem:interrupted-search} with
$M=n$ to this search. The lemma extends to a random finite stopping time
by applying it to finite budgets and taking monotone limits. This proves
the first bound in~\eqref{eq:coarse-lower}.

For supplied indices, use $k$ scalar coordinates, independent signs
$\sigma_1,\ldots,\sigma_k$, $A=nF_0/k$, and
$b^2=2A/(\delta\sqrt{k})$. Set
\[
 h(x)=\frac{L_0}{4}\|x\|^2,\qquad
 u_i(x)=-\sigma_i A\sin(x_i/b)\quad(i\leq k),\qquad
 u_i=0\quad(i>k).
\]
The same calculations, now using orthogonal coordinates, verify all
assumptions. Supply the exceptional indices and reveal a component's sign
on its first query. Conditional on the complete history, each unqueried
sign remains independent and uniform. The expected squared gradient in
its coordinate is at least $G/(2k)$, for any value of that coordinate
at the output. If $J$ distinct exceptional components are queried, then
\[
 \mathbb E\|\nabla F(\widehat x)\|^2
 \geq\frac{G}{2k}\bigl(k-\mathbb E J\bigr).
\]
Since the query count is at least $J$, this proves the second bound.
Both arguments average over fixed objective families; hence each expected
lower bound is attained by some fixed member.
\end{proof}

\subsection{Combining the bounds over the full accuracy range}
\label{app:phase-complexity}

\begin{proposition}[Explicit minimax lower bounds]\label{prop:explicit-minimax-lower}
In the setting of Theorem~\ref{thm:phase-complexity}, put
$D=1536\cdot1024$. Write $\mathfrak C_{\rm unknown}$ and
$\mathfrak C_{\rm supplied}$ for the two minimax costs. For
$0<\epsilon^2\leq\delta F_0\sqrt{k}/(4n)$,
\begin{equation}\label{eq:explicit-minimax-lower}
 \begin{aligned}
 \mathfrak C_{\rm supplied}&\geq2^{-23}\mathcal R,\\
 \mathfrak C_{\rm unknown}&\geq
 \max\left\{\frac{n+1}{4},\ 2^{-23}\mathcal R,\
       \frac{17}{8192D}\min\{nk,\mathcal P\}\right\}.
 \end{aligned}
\end{equation}
\end{proposition}

\begin{proof}[Proof of Proposition~\ref{prop:explicit-minimax-lower}
and Theorem~\ref{thm:phase-complexity}]
Write $C_*=1024$, $D=1536C_*$, and
\[
 G=\frac{\delta F_0\sqrt{k}}n,\qquad
 T=\frac G{\epsilon^2}=\frac{\mathcal R}{k},\qquad
 \vartheta=\frac{k}{T}=\frac{n\sqrt{k}\epsilon^2}{\delta F_0}
           =\frac{k^2}{\mathcal R}=\frac{nk}{\mathcal P}.
\]
The assumed accuracy range is $T\geq4$. Larger $\vartheta$ corresponds
to coarser accuracy.

\paragraph{Selecting an upper bound}
Choosing the smaller of the two upper bounds in
Section~\ref{sec:total-potential-budget} gives
$O(n+\min\{nk+\mathcal R,\mathcal P\})$ when $\Omega\leq2F_0$.
Both terms in the minimum are at least of order $n$ in the stated range.
Adding the cost of constructing the predictors therefore preserves this order.
Since $\mathcal P\geq\mathcal R$,
\[
 \min\{nk+\mathcal R,\mathcal P\}
 \leq\mathcal R+\min\{nk,\mathcal P\}
 \leq2\min\{nk+\mathcal R,\mathcal P\}.
\]

\paragraph{The accuracy lower bound}
If $T\geq D$, the parameter controlling chain length in
Theorem~\ref{thm:bounded-residual-lower} satisfies
\[
 A=\frac{T}{768\cdot512}=\frac{4T}{D}\geq4.
\]
That theorem gives the bound $2^{-23}\mathcal R$, even when the
exceptional indices are supplied. Its residual functions have total
oscillation at most $2nF_0$, as required. If $4\leq T<D$,
Proposition~\ref{prop:coarse-lower} gives a supplied-index lower bound
of at least $k/2$. Since $\mathcal R=kT<kD$ and
$1/(2D)>2^{-23}$, the same bound $2^{-23}\mathcal R$ follows.
Supplying additional information cannot increase the minimax cost, so this
lower bound also applies when the indices are unknown.

\paragraph{Coarse accuracy}
Proposition~\ref{prop:coarse-lower} gives
$\mathfrak C_{\rm unknown}\geq(n+1)/4$ throughout the stated range.
If $T<2D$, then $\mathcal P=nT<2Dn$, and hence
\[
 \mathfrak C_{\rm unknown}\geq\frac{\mathcal P}{8D}
 \geq\frac{17}{8192D}\min\{nk,\mathcal P\}.
\]
It remains to prove the discovery term for $T\geq2D$, or equivalently
$\vartheta\leq k/(2D)$.

\paragraph{Fine and intermediate accuracy}
If $\vartheta\leq2/D=1/(768C_*)$, the preceding discovery theorem
gives $\mathfrak C_{\rm unknown}\geq13nk/2048$. This implies the
claimed discovery term in~\eqref{eq:explicit-minimax-lower}.

\paragraph{Parallel chains at coarser accuracy}
Suppose $\vartheta>2/D$, and choose
\[
 m=\lceil D\vartheta\rceil,\qquad
 N=\lfloor n/m\rfloor,\qquad q=\lfloor k/m\rfloor.
\]
Here $m$ is the number of parallel chains, $N$ is the number of component
indices available to each chain, and $q$ is its length. The inequalities
$mN\leq n$ and $mq\leq k$ preserve the upper bounds on the numbers of
components and exceptions. Then
\[
 m\leq\tfrac32D\vartheta,\qquad
 m\leq k/2+1\leq k,\qquad N\geq q\geq1.
\]
Assign each chain a known group of $N$ component indices and a known
orthogonal coordinate block. Residual gradients at unused indices vanish identically.
In block $a$, independently draw a uniformly distributed orthonormal
$q$-frame $V_a$ and a uniform ordered list
$(\lambda_{a,1},\ldots,\lambda_{a,q})$ of distinct indices in its
component group. There are $k'=mq\leq k$ exceptions.

Let $d_1,\ldots,d_q$ be the links of Section~\ref{sec:siopt-chain}
with chain length $K=q$, so $c_q=\sum_jd_j$.
Set
\[
 \alpha=\frac{nF_0}{12k},\quad
 b^2=\frac{C_*\alpha}{\delta\sqrt{k}},\quad
 R_q=230\sqrt q,\quad
 \rho_a(y)=\frac{y}{\sqrt{1+\|y\|^2/R_q^2}},
\]
and define the fixed random family
\[
 \begin{aligned}
 h(x)&=\frac{\alpha\|x\|^2}{10nb^2},&
 u_{\lambda_{a,j}}(x)&=\alpha d_j(V_a^\top\rho_a(x_a/b)),\\
 u_i&=0\ \text{otherwise},& f_i&=h+u_i.
 \end{aligned}
\]
The objective still assigns weight one over $n$ to every component. In particular,
\[
 F(x)=\frac\alpha n\sum_{a=1}^m
 \left[c_q(V_a^\top\rho_a(x_a/b))+\frac{\|x_a/b\|^2}{10}\right].
\]

\paragraph{Verifying the assumptions after rounding}
The elementary link bounds from Lemma~\ref{lem:siopt-chain}
are $\|\nabla d_j\|\leq23$, $\operatorname{osc}(d_j)\leq24$, and
$\|\nabla^2d_j\|_{\rm op}\leq410$; each coordinate occurs in at most
two links. Within each block, the chain-rule calculation used above bounds
the Euclidean norm of the concatenated Hessian-vector products by
\[
 \left(\sqrt2\,410+\frac{6\cdot23\sqrt q}{R_q}\right)\|v_a\|
 <C_*\|v_a\|.
\]
The coordinate blocks are orthogonal, so summing $\|v_a\|^2$ across
blocks gives $\|v\|^2$, without a factor $m$. The chosen scale therefore
gives the required residual smoothness bound $k\delta^2\|x-y\|^2$
with the supplied $k$, although only $k'$ components
are exceptional. The gap is at most
$12mq\alpha/n=F_0k'/k\leq F_0$. The objective Hessian is block diagonal
and has norm at most $512\alpha/(nb^2)\leq L_0$.
Moreover
\[
 \|r_i(x)\|^2\leq\frac{529\alpha^2}{b^2}
   =\frac{529}{12C_*}\frac{\delta nF_0}{\sqrt{k}},
 \qquad \operatorname{osc}(u_i)\leq24\alpha=\frac{2nF_0}{k}.
\]
Consequently $\sum_i\operatorname{osc}(u_i)\leq2nF_0k'/k\leq2nF_0$.
All bounds hold on the whole space. The algorithm may know the partition,
the unused indices, and the quadratic reference; this information does not
change the argument.

\paragraph{Agreement for adaptively chosen queries}
In each block, keep track of the successive links already revealed. The auxiliary oracle returns
exact gradients on revealed indices and on the next index in each incomplete
block, revealing its index and frame column on that query; for all other
indices not yet revealed, it returns only $\nabla h(x)$. Define this process on all
outcomes, and run the given algorithm while withholding the additional
revelations. We next give an event on which its gradient responses agree
with those of the actual oracle.

Let $Q=\lfloor mNq/64\rfloor$. Give every coordinate block dimension
\[
 d_0\geq q+1+8R_q^2\log\bigl(128k(Q+1)\bigr).
\]
Conditional on the auxiliary history, the unrevealed columns form a
uniformly distributed orthonormal frame in the complement of the columns
already revealed in their block. Under this conditioning, the frame is
independent of the ranks of the unknown component indices. We do not condition this distribution
on successful projection tests. At each query and at the output, test the
projections onto the unrevealed columns. There are at most $k(Q+1)$ tests.
For each test, concentration for a uniformly distributed unit vector bounds
the probability of projection magnitude at least $1/2$ by
$2\exp[-(d_0-q-1)/(8R_q^2)]$. The probability that any test fails is at most $1/64$. If all tests
succeed, the coordinate preceding each later link is unrevealed and has
magnitude below $1/2$. That link and its gradient vanish near the query.
The actual and auxiliary processes then have identical queries, responses,
and outputs. This includes
an arbitrary unqueried output and permits unbounded adaptive query points.

\paragraph{Searching across blocks in any adaptive order}
Within a block, each successive label is conditionally uniform among
its remaining indices. The all-failure continuation in
Lemma~\ref{lem:interrupted-search} is independent of that label: replies
from earlier links use already revealed directions, and other blocks
have independent assignments and frames. The lemma therefore applies to
the auxiliary process, including arbitrary switches between blocks.
With $C_Q$ the number of completed blocks, it gives
\[
 \mathbb E C_Q\leq\frac{4Q}{Nq}\leq\frac m{16},\qquad
 \Pr\{C_Q\geq m/2\}\leq\frac18.
\]
This argument holds for every $q\geq1$. In particular, it does not
require a union bound over subsets of completed blocks.

\paragraph{The gradient bound and random stopping times}
With probability at least $55/64$, the actual and auxiliary processes
agree and at least $m/2$ blocks have unrevealed links. In each such block,
the gradient bound of Lemma~\ref{lem:siopt-chain} gives
$\|\nabla_a F(\widehat x)\|>\alpha/(2nb)$.
Orthogonality yields
\[
 \|\nabla F(\widehat x)\|^2>
 \frac{m\alpha^2}{8n^2b^2}
 =\frac{m\delta F_0}{96C_*n\sqrt{k}}
 =16\epsilon_m^2,\qquad
 \epsilon_m^2=\frac{m\delta F_0}{Dn\sqrt{k}}\geq\epsilon^2.
\]
Write $N_{\rm paid}$ for the
original algorithm's total number of component-gradient queries.
Let $\widehat x_Q$ denote the
output of an almost-surely terminating algorithm run until return or
immediately before a requested $(Q+1)$st gradient, with zero returned in
the latter case. It agrees with the original output
$\widehat x$ on $\{N_{\rm paid}\leq Q\}$. The assumed bound on the
expected squared gradient norm at $\widehat x$ and Markov's inequality give
\[
 \frac{55}{64}\leq\Pr\{\|\nabla F(\widehat x_Q)\|^2\geq16\epsilon^2\}
 \leq\Pr\{N_{\rm paid}>Q\}+\frac1{16}.
\]
Hence $\Pr\{N_{\rm paid}>Q\}\geq51/64$. Using the integer query
count avoids a restriction on the size of $mNq$:
\[
 \mathbb E N_{\rm paid}
 \geq(Q+1)\Pr\{N_{\rm paid}>Q\}
 \geq\frac{51mNq}{4096}.
\]
Since $mN\geq n/2$, $mq\geq k/2$, and
$m\leq3D\vartheta/2$, we obtain
\[
 \mathbb E N_{\rm paid}\geq\frac{51nk}{16384m}
 \geq\frac{17}{8192D}\frac{nk}{\vartheta}
 =\frac{17}{8192D}\mathcal P.
\]
Averaging over the objective family gives a fixed member attaining this
cost. Together with the coarse-accuracy and single-chain cases, this proves
\eqref{eq:explicit-minimax-lower}. The maximum of its accuracy and discovery
terms bounds a positive constant times
$\mathcal R+\min\{nk,\mathcal P\}$, completing the lower bound in the main
theorem. Different families may attain its separate terms; all belong to
the same objective class.

\paragraph{Known exceptional indices and the accuracy ranges}
When the exceptional indices are supplied, the residual estimator costs
$O(k+\mathcal R)$. It checks the exact gradient using only those indices,
since all omitted residual gradients vanish identically. Since
$\mathcal R\geq4k$, this is $O(\mathcal R)$; the accuracy lower bound
above gives a matching lower bound. If the query convention also requires querying all $n$ components,
the total cost is $\Theta(n+\mathcal R)$. Otherwise, this additional
cost is absent. Finally,
$\mathcal P=(n/k)\mathcal R$ crosses $nk$ at $\mathcal R=k^2$, and
$\mathcal R$ crosses $nk$ at $\mathcal R=nk$, which gives the stated accuracy ranges.
\end{proof}
\section{Constructing predictors from known loss functions}
\label{sec:predictor-acquisition}
We construct predictors as mixtures of known loss gradients. This construction
is related to mixture models for client distributions~\cite{marfoq2021mixtures},
estimation from reference gradients~\cite{wang2021imbalance}, and the use of
moment equations to identify label shift~\cite{lipton2018labelshift}.
The mixture assumption must hold at every point; the queries determine its
coefficients.

\begin{proposition}[Constructing predictors under a mixture assumption]
\label{prop:acquired-mixture-predictor}
Let $\mathsf r$ be a known differentiable regularizer and $\psi_1,\ldots,\psi_m$ known
differentiable losses, where $m\geq2$. Write
$\mathsf G(x)=[\nabla\psi_1(x),\ldots,\nabla\psi_m(x)]$ and
$\Delta_m=\{\theta\geq0:\mathbf1_m^\top\theta=1\}$.
Suppose that, outside an unknown set of at most $k$ components,
\[
 \nabla f_i(x)=\nabla\mathsf r(x)+\mathsf G(x)\theta_i+e_i(x),
 \qquad \theta_i\in\Delta_m,\qquad \sup_x\|e_i(x)\|\leq\zeta.
\]
For specified points $z_1,\ldots,z_q$ and an orthonormal basis
$U\in\mathbb R^{m\times(m-1)}$ of $\mathbf1_m^\perp$, assume
\[
 D=q^{-1/2}\operatorname{stack}_{\ell=1}^q\mathsf G(z_\ell),
 \qquad \kappa=\sigma_{\min}(DU)>0,
 \qquad H_D=\sup_x\|\mathsf G(x)U\|_{\rm op}<\infty.
\]
Here $\operatorname{stack}$ denotes vertical concatenation.
Query every component at these points and set
\[
 y_i=q^{-1/2}\operatorname{stack}_{\ell=1}^q
       [\nabla f_i(z_\ell)-\nabla\mathsf r(z_\ell)],\qquad
 \widehat\theta_i=\mathop{\rm argmin}_{\theta\in\Delta_m}
                      \|D\theta-y_i\|^2.
\]
Set $p_i=\mathsf r+\sum_j\widehat\theta_{ij}\psi_j$, so that
$\nabla p_i=\nabla\mathsf r+\mathsf G\widehat\theta_i$. These predictors satisfy
\[
 \|\widehat\theta_i-\theta_i\|\leq\zeta/\kappa,
 \qquad \sup_x\|\nabla f_i(x)-\nabla p_i(x)\|
       \leq\zeta(1+H_D/\kappa)
\]
on every component satisfying the mixture assumption. Constructing these
predictors costs $C_{\rm cal}=nq$ queries. In particular, $\zeta=0$ gives exact predictors
outside at most $k$ components.
\end{proposition}
\begin{proof}
Write $\theta=\mathbf1_m/m+Uc$ and let $C$ denote the simplex in these
coordinates. The least-squares fit projects $y_i-D\mathbf1_m/m$ onto the
closed convex set $DUC$. The stacked observation error has norm at most
$\zeta$. Projection does not increase Euclidean distances, so comparison
with the true coefficient gives $\|DU(\widehat c_i-c_i)\|\leq\zeta$.
The smallest singular value of $DU$ bounds the coefficient error, and
$e_i(x)+\mathsf G(x)U(c_i-\widehat c_i)$ is the prediction error.
\end{proof}

The condition $\kappa>0$ requires $qd\geq m-1$, but this dimension
inequality alone does not ensure that the coefficients can be recovered.
With $\bar\theta=n^{-1}\sum_i\widehat\theta_i$, the average predictor is
$\nabla h=\nabla\mathsf r+\mathsf G\bar\theta$. Evaluating it uses $m$
known loss gradients, which are uncharged in our query model. A component
query may still require a full pass through that component's data.
The $nq$ queries used to construct the predictors do not change the order
of the main complexity bound when $q$ is bounded independently of the
problem parameters.

\paragraph{Bounds on residual oscillation and smoothness}
Suppose $f_i=\mathsf r+L_i$ for every component, with
$0\leq L_i,\psi_j\leq M$ and $L_\ell$-Lipschitz loss gradients. Then
\[
 u_i=L_i-\sum_j\widehat\theta_{ij}\psi_j,\qquad
 \operatorname{osc}(u_i)\leq2M,\qquad
 \operatorname{Lip}(\nabla u_i)\leq2L_\ell.
\]
For $\zeta=0$, only the exceptional components have nonconstant residual
functions, giving $\Omega=2Mk/n$ and $\delta=2L_\ell$. For $\zeta>0$,
the uniform prediction bound remains valid, but more than $k$ residual
gradients may be nonzero.

For example, take $d\geq m$, $\psi_j(x)=\tfrac12(\tanh x_j-1)^2$, and
$z_1=0$. Then $\mathsf G(0)=-[e_1,\ldots,e_m]$, $\kappa=1$, and
$H_D\leq2$. If every component loss, apart from its regularizer, is an
average of
$\tfrac12(\tanh(a^\top x)-y)^2$ with $\|a\|\leq1$ and
$y\in\{-1,1\}$, valid bounds are $M=2$ and $L_\ell=5$.
With $\mathsf r=\lambda\|x\|^2/2$, $\lambda\geq0$, and $x_0=0$,
we may take $L=\lambda+5$, $F_0=1/2$, and $\delta=10$.
If the exact mixture assumption holds outside at most $k\leq n/4$
components, then
$\sum_i\operatorname{osc}(u_i)\leq4k\leq2nF_0$.

\begin{proposition}[Identifying analytic residuals with one query per component]
\label{prop:analytic-support-probe}
Suppose the residual gradients $r_i$ are real analytic on $\mathbb R^d$.
Query every component at a common independent point drawn from a distribution
with a Lebesgue density. Testing $r_i(x)=0$ identifies $\mathcal B$ almost
surely. Under the other assumptions of Theorem~\ref{thm:phase-complexity},
optimization then costs $O(C_{\rm cal}+n+\mathcal R)$ queries.
\end{proposition}
\begin{proof}
Each residual gradient that is not identically zero has a scalar coordinate
that is a nonzero analytic function. Its zero set has Lebesgue measure
zero~\cite{mityagin2020zero}. The finite union of these zero sets also has
measure zero. Proposition~\ref{prop:tracking-population} then applies to
the identified exceptional set, including the final exact gradient check.
\end{proof}
This result uses exact arithmetic. A positive detection threshold or noisy
response does not give the same guarantee. The preceding analytic example
therefore illustrates optimization after the exceptional components have
been identified. The mixture assumption must still hold globally; these
finitely many queries do not establish it.

\paragraph{A smooth loss with flat regions}
\label{ex:nonanalytic-margin}
The same predictor construction permits smooth losses that are constant
on open sets.
Define
\[
 \ell(s)=
 \begin{cases}
 \exp\!\bigl(1-(1-2s)^{-2}\bigr),&s<1/2,\\
 0,&s\geq1/2.
 \end{cases}
\]
This $C^\infty$ loss takes values in $[0,e]$ and satisfies
$\ell(0)=1$, $\ell'(0)=-4$,
$\|\ell'\|_\infty\leq\sqrt{54/e}$, and $\|\ell''\|_\infty<33$.
For $\psi_j(x)=\ell(x_j)$ and $d\geq m$, the origin gives
$\kappa=4$ and $H_D\leq\sqrt{54/e}$.
Let each $L_i$ be an average of $\ell(ya^\top x)$ with $\|a\|\leq1$
and $y\in\{-1,1\}$. Assume that each $L_i$ is an exact mixture of the
known losses, except on at most $k$ components. With $f_i=\lambda\|x\|^2/2+L_i$, $\lambda\geq0$,
valid bounds at $x_0=0$ are
\[
 C_{\rm cal}=n,\qquad L=\lambda+33,\qquad \delta=66,
 \qquad F_0=1,\qquad \sum_i\operatorname{osc}(u_i)\leq2ek.
\]
Thus $k\leq n/e$ gives the required bound on the sum of residual oscillations.

For an exceptional component, take $m=d=2$ and
\[
 L_i(x)=\ell\!\left(\frac{x_1+x_2}{2}\right),\qquad
 u_i(x)=L_i(x)-\frac{\ell(x_1)+\ell(x_2)}2.
\]
Fitting the coefficients from a query at the origin gives $(1/2,1/2)$
and $\nabla u_i(0)=0$. Although $\nabla u_i(2,0)=2e_2$, both $u_i$ and
its gradient vanish on $\{x_1>1/2,x_2>1/2\}$. A random query point that
has positive probability of lying in this open set can fail to detect the
exception. Thus a smooth mixture model can have residual gradients for
which the identification guarantee for analytic functions does not apply.
\section{Numerical experiments}
\label{sec:experiments}
We examine delayed identification of exceptional components and optimization
once their indices have been found. The first study below tests the adaptive
working-set method against baselines that use all support information they
acquire. A separate regression study changes the representation of a fixed
objective after its exceptional indices are identified. We distinguish
component queries, optimization steps, and private gradient work.

\subsection{Delayed discovery with nonzero estimator error}
\label{sec:delayed-experiments}
We first examine whether adapting the batch size helps when exceptions
become visible at different stages of optimization. Every method starts
without the exceptional indices, except for a separately reported
supplied-support comparison. The component residuals vanish at the initial
point, so an initial scan cannot identify them.

\paragraph{Objective and valid bounds}
Let $x\in\mathbb R$, $x_0=0$, and $h(x)=(x-1)^2/2$. Define
\[
 q(z)=\begin{cases}
 0,&z\leq0,\\
 10z^3-15z^4+6z^5,&0<z<1,\\
 1,&z\geq1.
 \end{cases}
\]
The $k$ exceptional components form $J$ equal groups, with indices
chosen uniformly without replacement from $[n]$. For group
$j\in\{0,\ldots,J-1\}$, put
\[
 t_j=j+\tfrac12,\qquad T_j=j+J+3,\qquad
 g_j(x)=q(4(x-t_j))-q(4(x-T_j)).
\]
Within each group, components occur in pairs. Draw a weight in $[0.9,1.1]$
and a phase in $[0,2\pi]$ uniformly and independently for each pair.
Give its two members equal weights and phases and opposite signs
$\sigma_i\in\{-1,1\}$. Normalize the weights $w_i$ to sum to one in
each group. For a component in group $j$, set
\begin{equation}\label{eq:delayed-residual-family}
 r_i(x)=-nw_i g_j(x)\bigl(1+\beta\sigma_i\sin(x+\phi_i)\bigr),
 \qquad \beta=0.3,
\end{equation}
and let $u_i$ be its integral, equal to zero below $t_j$.
Other residual functions vanish. All predictors are $p_i=h$.
The paired sine terms cancel in the group sum. Consequently
\[
 F'(x)=x-1-\sum_{j=0}^{J-1}g_j(x),
\]
which is independent of $n$, the hidden indices, the weights, and the
phases. Exceptions in group $j$ are zero for $x\leq t_j$ and become
visible near the partial stationary point $j+1$.
The later decrease of $g_j$ to zero makes each residual potential bounded.

The following global bounds are supplied to the methods:
\begin{equation}\label{eq:delayed-family-bounds}
 \begin{aligned}
 L&=8.5,&
 \delta&=\frac nk\frac{11}{9}
       \left[7.5(1+\beta)\sqrt J+\beta J\right],\\
 \Omega&=(1+\beta)J(J+2.5),& F_0&=\tfrac12+\Omega.
 \end{aligned}
\end{equation}
Indeed, $0\leq q'\leq15/8$, the transition intervals are disjoint, and
$w_i\leq(11/9)J/k$. The derivative of the stacked residual vector has
norm at most
$n(11/9)[7.5(1+\beta)\sqrt{J/k}+\beta J/\sqrt k]$,
which gives the stated $\delta$. Each residual function is nonpositive and
monotone, with oscillation at most
$(1+\beta)nw_i(T_j-t_j)$. Summing gives $\Omega$, and $F\geq-\Omega$
gives $F_0$. Thus $\Omega\leq2F_0$.
Every full-complement check adds complete groups because the factors in
parentheses in~\eqref{eq:delayed-residual-family} lie in $[0.7,1.3]$.
Each working set is therefore a union of complete groups. The sine terms
cancel in every partial objective the working-set methods can form, giving
the global bound $\overline L=8.5$ for these objectives. Arbitrary subsets
would instead require the larger valid bound
$1+7.5(1+\beta)+\beta J$.

\paragraph{Methods and information}
The adaptive method A is Theorem~\ref{thm:adaptive-working-set} with
$\eta=1/(2\overline L)=1/17$ and
$p=\min\{1,2\eta\delta\sqrt{k}/n\}$. Its expected batch size is
$\eta\delta\sqrt{k}|S|/n$.
F uses the same working set, recomputation probability, stepsize, and
exclusive choice between recomputation and sampled differences, with the
fixed expected batch size $\eta\delta\sqrt{k}k/n$.
GD recomputes the working-set gradient at every step. The final comparison
uses its standard full step $\eta=1/\overline L=1/8.5$, justified by the
exact-gradient descent bound in Section~\ref{sec:adaptive-working-set}.
These three methods use the same full-complement discovery rule, the
threshold $\epsilon$, reuse of already queried gradients, and immediate
certification once $|S|=k$.
P uses population sampling from Proposition~\ref{prop:tracking-population},
with $\eta=1/17$. It records every observed nonzero residual,
including both endpoints of sampled differences, and restricts sampling to
the support at the next exact recomputation once $k$ exceptions are known.
The supplied-support method O is given the true indices before its first
query; it is excluded from comparisons between methods with equal
information. We also rerun A with the more conservative choice
$\eta=1/[2(L+\mathcal H)]$ to measure the effect of the supplied
partial smoothness bound.

\paragraph{Study sequence and measurement}
We conducted five sequential studies, fixing each source, protocol, and
set of evaluation seeds before running that study. The first used constant
residual gradients between discoveries: eight settings, eight methods,
and twelve seeds gave 768 outcomes. The second introduced the paired
sine terms and retained the conservative tracking stepsize: five settings,
five methods, and twelve new seeds gave 300 outcomes. Its estimator errors
were nonzero, and its adaptive method used more queries than GD in all
five settings. The third used the supplied partial smoothness bound:
six methods on another twelve seeds gave 360 outcomes. That comparison
used $\eta=1/17$ for both tracking and GD. We then added GD with its
larger valid step $1/8.5$ on the same sixty instances. This changed the
comparison: full-step GD used fewer queries than A in every instance.
All preceding protocols and results are retained.

The first three studies use $k=32$. The third has
$n\in\{256,2048\}$ crossed with $J\in\{4,8\}$, together with
$n=k=32$, $J=4$. Its development seeds are $51410,51411$ and its
evaluation seeds are $51600,\ldots,51611$; the GD addendum uses the same
evaluation seeds. The earlier studies have disjoint seed sets.

The remaining study tests the predicted dependence on support size.
At fixed $J$ and $n/k=8$, the residual scale $\delta\sqrt{k}/n$ and
recomputation probability are proportional to $k^{-1/2}$, while
$\overline L=8.5$ stays fixed. Before running new evaluation seeds, we
fixed $k\in\{128,512\}$, $n=8k$, $J\in\{4,8\}$, and the four methods
A, F, full-step GD, and P. Twelve new seeds, $51700,\ldots,51711$,
give 192 outcomes. No parameter grid was tuned on those outcomes.
For $J=8$, the recomputation probabilities at $k=32,128,512$ are
approximately $0.762$, $0.381$, and $0.191$. Their expected update costs
per visit, as fractions of exact recomputation on $S$, are
$(2-p)p$, approximately $0.943$, $0.617$, and $0.345$.
The total cost also depends on the steps and discovery checks.

Across all studies, $\epsilon=10^{-4}$ and the optimizer's random seed
is the instance seed plus 9871. Each run allows $10^5$ steps and $2000n$
queries. No method uses momentum or line search. All 1680 outcomes pass
a paid full-gradient certificate; none reaches a limit. Independent
evaluation of the explicit formula for $F'$ verifies every returned point.
The code, protocols, seeds, and individual outcomes accompany the manuscript.

Every component query is charged, including initial queries, both endpoints
of a sampled difference, repeated indices, and final certificates. We
separately count evaluations of residual gradients on exceptional
components as private work. This count assumes no caching and does not
measure execution time. Clipped polynomial evaluations give exact zeros
in the designed inactive regions; the implementation tests a response
against zero in double precision. Estimator-error calculations are performed
only after a run and do not inform optimization.

\paragraph{Results and limits}
Table~\ref{tab:delayed-public-smoothness} reports all nine settings with
full-step GD. Adaptive batching reduces mean queries relative to fixed
batching in every setting. Its comparison with GD depends on the residual
scale. GD uses fewer queries in all five $k=32$ settings and at
$k=128$, $J=8$, where A costs 23.9\% more. A uses 3.4\% fewer queries
at $k=128$, $J=4$, and 41.7\% and 23.0\% fewer at $k=512$, $J=4$ and
$J=8$, respectively. Each comparison holds in all twelve paired runs of
its setting. P costs least in the dense case. These results distinguish
the benefit of adapting the batch from the separate choice between
stochastic estimation and exact gradients.

\begin{table}[htbp]
\caption{Paid queries with full-step GD: mean (empirical standard deviation),
rounded to the nearest query, over twelve paired seeds per setting.
A adapts the batch size, F fixes it, GD uses exact working-set gradients,
and P samples the population and then the identified support.
A, F, and P use $\eta=1/17$; GD uses $\eta=1/8.5$.
The five $k=32$ rows combine the third study with its paired GD addendum;
the four larger-support rows use new seeds in the final study.
All returned points have gradient norm at most $10^{-4}$.}
\label{tab:delayed-public-smoothness}
\centering\small
\begin{tabular}{@{}rrr rrrr@{}}
\toprule
$n$ & $k$ & $J$ & A & F & GD & P\\
\midrule
32 & 32 & 4 & 9528 (150) & 11413 (129) & 6000 (0) & 5906 (133)\\
256 & 32 & 4 & 10437 (146) & 12338 (133) & 6896 (0) & 19189 (714)\\
256 & 32 & 8 & 22667 (83) & 25734 (161) & 12592 (0) & 43787 (319)\\
2048 & 32 & 4 & 17605 (158) & 19496 (153) & 14064 (0) & 127674 (3618)\\
2048 & 32 & 8 & 36976 (71) & 40058 (186) & 26928 (0) & 314291 (3364)\\
1024 & 128 & 4 & 26644 (461) & 32407 (512) & 27584 (0) & 46450 (3036)\\
1024 & 128 & 8 & 62414 (653) & 78615 (562) & 50368 (0) & 117001 (4352)\\
4096 & 512 & 4 & 64302 (2291) & 78980 (2478) & 110336 (0) & 101318 (9540)\\
4096 & 512 & 8 & 155118 (4114) & 197490 (3432) & 201472 (0) & 262789 (13656)\\
\bottomrule
\end{tabular}
\end{table}

Every run of A discovers one complete group at each failed check, so the
observed number of discovery stages is $D=J$. Its complementary scans use
exactly $Jn-k(J-1)/2$ queries. Individual sampled differences and estimator
errors are nonzero. The exact recomputations and terminal checks remain
essential even on this controlled family. The preceding studies retain
both the common-stepsize results and the unfavorable results for the
conservative tracking stepsize.

\begin{figure}[htbp]
\centering
\includegraphics[width=\linewidth]{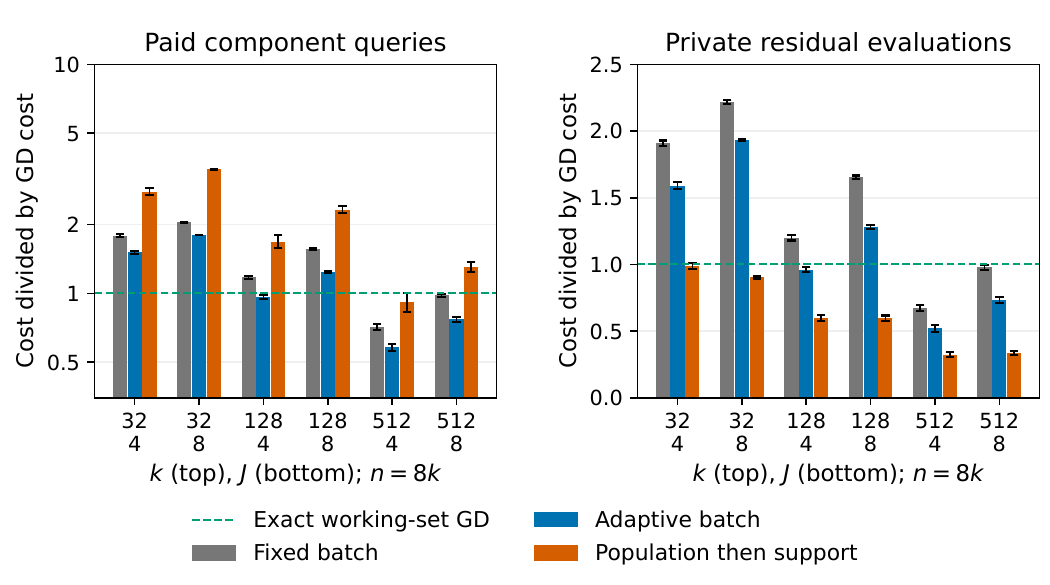}
\caption{Queries and private residual work divided by the cost of
full-step working-set GD, for all settings with $n/k=8$.
Bars show arithmetic means over twelve paired seeds; error bars show one
empirical standard deviation divided by the deterministic GD cost.
Ratios below one indicate lower cost. A and F use $\eta=1/17$; GD uses
$\eta=1/8.5$. Population sampling may use fewer private evaluations even
when it makes more component requests. The table also includes the dense
case and the two settings with $n/k=64$.}
\label{fig:delayed-resources}
\end{figure}

Figure~\ref{fig:delayed-resources} compares query and private-work costs
at fixed $n/k=8$. At $n=4096$, $k=512$, $J=8$, A averages 155,118.3
queries and 126,446.3 private residual evaluations, compared with
262,789.0 queries and 58,041.7 private evaluations for P.
A saves 41.0\% of the requests while performing 2.18 times as much
private residual work. It averages 1409.8 steps, compared with 324.2 for
P and 666 for full-step GD. Against GD, A reduces queries by 23.0\% and
private residual evaluations by 26.8\%. The supplied-support comparison
remains separate: for example, at $n=2048$, $k=32$, $J=8$, O averages
9728.3 queries when its exceptional indices are supplied initially.

These examples satisfy the main theorem's accuracy condition, but do not
approximate its worst-case high-dimensional constructions. They are
one-dimensional, use exact gradient access and a supplied partial smoothness
bound, and have a partial residual sum that is constant between activation
intervals. The sine terms test nonzero component differences and estimator
error while cancelling in the full objective. The results establish a
benefit for the specified discovery mechanism and parameter choice; they
do not establish performance on general regression data, noisy residuals,
or distributed implementations.

\subsection{Regression objective and protocol}
\begin{proposition}[Changing the components without changing the objective]
\label{prop:fixed-regression-objective}
Let $\|a_j\|\leq1$, $\phi_j(x)=\tfrac12\tanh^2(a_j^\top x-b_j)$,
and $h(x)=\tfrac\mu2\|x-c\|^2$, where $\mu>1/3$.
For $2\leq k\leq n$ with $k\mid M$, choose $A\subseteq[n]$, $|A|=k$,
and partition $M$ examples into equal blocks $(B_i)_{i\in A}$. Define
\[
 \varphi_i=\frac{k}{M}\sum_{j\in B_i}\phi_j,\qquad
 u_i=\frac nk\varphi_i\quad(i\in A),\qquad
 u_i=0\quad(i\notin A).
\]
The true exceptional set satisfies $\mathcal B\subseteq A$.
With $p_i=h$, the objective $F=h+M^{-1}\sum_j\phi_j$ is independent of
$n$, $k$, and the partition. Valid bounds are
\[
 L=\mu+1,\quad \delta=n/k,\quad \Omega=1/2,\quad
 F_0=\tfrac\mu2\|x_0-c\|^2+1/2.
\]
The objective is $(\mu-1/3)$-strongly convex, and the partial objective
for every working set is $(\mu+1)$-smooth.
\end{proposition}
\begin{proof}
The scalar loss has range $[0,1/2)$ and second derivative
$(1-t^2)(1-3t^2)\in[-1/3,1]$, where $t=\tanh z$.
Hence $-I/3\preceq\nabla^2\phi_j\preceq I$; averaging proves the
objective and smoothness statements. Each residual function has oscillation at most
$n/(2k)$ and gradient Lipschitz constant at most $n/k$.
Summing these bounds proves the residual assumptions. Nonnegativity of $F$
gives the bound on the initial objective gap.
\end{proof}

We use $M=1024$, $d=4$, $\mu=0.5$, $c=x_0=0$, and $\epsilon=10^{-4}$.
Each feature is the vector $e_1+0.25z_j$, $z_j\sim N(0,I)$, normalized
to have unit norm.
Half the targets have mean $0.05$ and half have mean $0.7$, with independent
$N(0,0.015^2)$ perturbations. Thus $L=1.5$ and $F_0=\Omega=0.5$.
Outputs that pass the gradient check lie within $6\times10^{-4}$ of
the unique minimizer
and have objective gap at most $3\times10^{-8}$.

A clustered partition sorts the examples by target value and divides them
into consecutive blocks. A mixed partition assigns the sorted examples to
blocks in cyclic order.
We vary $n\in\{128,512,2048,8192\}$ with $k=64$ and clustered blocks,
and compare both partitions for $k\in\{16,64,256\}$ at $n=2048$.
One configuration occurs in both comparisons, giving nine distinct settings.
For each data seed, all settings use the same examples and objective.
Increasing $n$ adds components whose residual gradients vanish and rescales
the nonzero residual functions. The bound $\delta=n/k$ therefore changes
with the number of components.

\paragraph{Methods and stopping}
This regression comparison fixes each method's sampling domain to separate
the cost of estimating a residual sum from the cost of querying zero
residuals. Its population methods continue to sample all components even
after the initial scan identifies the support. The comparison in
Section~\ref{sec:delayed-experiments} allows population sampling to use
that information and switch to the identified support.
Every method first queries all $n$ components at $x_0$.
P-GD applies gradient descent using all components at each step; S-GD
applies gradient descent using the identified exceptions. S-Track and
P-Track use the estimator in Proposition~\ref{prop:tracking-population},
sampling from the $k$ identified exceptions and all $n$ components,
respectively. Between exact recomputations, they update an
estimate of the residual gradient sum using sampled differences.
These four methods use $\eta=1/(2L)=1/3$.
The fifth method, P-One, uses the same population estimator with at most one
sampled difference per step and $\eta=\min\{1/(2L),1/(\ell_A\sqrt n)\}$.
For a sampling set of size $m$, the estimator parameters are
\[
 \ell_A=\sqrt{mk}\delta/n,\qquad
 p_A=\min\{1,2\eta\ell_A/\sqrt m\},\qquad
 \tau_A=\eta\ell_A\sqrt m.
\]
The samples and independent Bernoulli variables that decide which differences
to evaluate follow Proposition~\ref{prop:tracking-population}.
S-GD and S-Track use the identified set only if the initial scan finds $k$
nonzero residual gradients, reaching the supplied upper bound on their number.
Otherwise, they use all components. The initial queries identify all $k$
exceptions in every run. The threshold for detecting a nonzero residual
gradient is $10^{-13}$. A method stops only after an exact gradient check
finds a norm at most $\epsilon$, initially or after a later exact
recomputation. No method uses momentum or line search. We count evaluations
at both points of every sampled difference, including repeated component
indices. The final check uses the gradient already computed.

\paragraph{Repetitions and resources}
We evaluate all five methods in all nine settings with the same eight data
seeds, $43100,\ldots,43107$. The protocol was fixed using the separate
seeds $42100,42101$. NumPy 2.0.2's default generator draws features followed
by target noise. For data seed $s$, another generator with seed $s+1009k$
randomly permutes $\max\{128,2k\}$ component positions and assigns the
data blocks to the first $k$ positions in that permutation. The optimizer
knows neither these positions nor the restriction used to choose them.
Random choices within each method use seed $s+781$. Each run allows
200,000 steps and $2000n$ queries. All 360 runs pass the stopping gradient
check. Full gradients computed after each run independently verify the
outputs and do not inform any iteration. We report arithmetic means and
observed ranges across seeds; these summaries are descriptive.

A query to a component containing a data block requires $M/k$ example
gradient evaluations. A query to a component containing no data block evaluates only the known
gradient $\nabla h$. We count repeated evaluations separately even
when the simulator combines repeated requests to the same index. This
count describes the computations required when each request is evaluated
separately; it does not measure simulator runtime, communication, or evaluations of
the known functions.

\subsection{Query savings after identifying the exceptions}
Table~\ref{tab:controlled-main} includes the common initial scan.
After these initial 2048 queries, the mean additional counts are 43,264
for P-GD, 1352 for S-GD, and 603.5 for S-Track. The first reduction avoids queries to zero
residual gradients. The second avoids repeatedly summing all exceptional
residual gradients. S-Track uses 22.0\% fewer total queries than S-GD and
86.7\% fewer than P-Track, with more steps than S-GD. Its mean count of example-gradient evaluations is 52.9\% lower than
that of S-GD.

\begin{table}[t]
\caption{Clustered blocks, $n=2048$, $k=64$: means [minimum--maximum]
over eight seeds. Queries include the initial scan and final gradient
check. The last column gives the mean count of example-gradient evaluations,
counting repeated evaluations separately.}
\label{tab:controlled-main}
\centering\small
\begin{tabular}{lrrr}
\toprule
Method & Component queries & Steps & Example gradients\\
\midrule
P-GD & 45{,}312.0 [43{,}008--47{,}104] & 21.1 [20--22] & 22{,}656.0 \\
S-GD & 3{,}400.0 [3{,}328--3{,}456] & 21.1 [20--22] & 22{,}656.0 \\
S-Track & 2{,}651.5 [2{,}490--2{,}844] & 66.4 [39--119] & 10{,}680.0 \\
P-Track & 19{,}866.5 [13{,}822--26{,}124] & 53.5 [25--86] & 9{,}800.0 \\
P-One & 25{,}380.2 [16{,}238--41{,}902] & 6{,}551.1 [3{,}003--15{,}835] & 12{,}632.0 \\
\bottomrule
\end{tabular}
\end{table}

For each data seed, adding components whose residual gradients vanish
leaves the iterates and query count of S-GD unchanged after the initial
scan. The same statement holds separately for S-Track. The mean total counts are $n+1352$ for S-GD and
$n+603.5$ for S-Track.
The relative saving therefore decreases from 50.6\% at $n=128$ to 7.8\%
at $n=8192$, where the initial scan accounts for 93.1\% of S-Track's cost.

Both estimators recompute their complete residual gradient sums with
probability $p=\min\{1,2\eta/\sqrt k\}$. Conditional on making no exact recomputation, their expected numbers of
sampled differences are $\tau_s=\eta\sqrt k$ and
$\tau_p=\eta n/\sqrt k$. Thus the expected query counts per step satisfy
$c_s=kp+2(1-p)\tau_s$ and $c_p=(n/k)c_s$.
Each exact recomputation evaluates all $M$ examples. A component sampled uniformly from all $n$ components contains a data
block with probability $k/n$. The two methods therefore have the same expected count of example-gradient
evaluations per step. Their total counts can differ
because they stop at different times.

\subsection{Effects of partitioning and batching}
Table~\ref{tab:controlled-partition} reports all nine settings.
Both gradient descent methods take 21.125 steps on average in every setting,
giving mean query
costs $22.125n$ for P-GD and $n+21.125k$ for S-GD.
Across the six partition settings, S-Track reduces mean total queries
relative to S-GD by 6.9--46.6\%.

To compare gradient changes across components at a fixed pair of points,
set $x=0$,
$y=0.2e_1$, $D_i=\nabla\varphi_i(y)-\nabla\varphi_i(x)$, and
$V=k^{-1}\sum_i\|D_i-k^{-1}\sum_jD_j\|^2$.
Mixing reduces mean $V$ by factors of 9.5--151.7, but at $k=256$ raises
S-Track's mean queries from 3985.0 to 4126.5.
The estimator randomly omits some candidate differences. Its variance
therefore also depends on the squared norms of the individual gradient
changes, before their mean is subtracted. Its stopping time depends on
when exact recomputations occur. A smaller value of $V$ at this pair of
points consequently need not give a smaller total query count.

\begin{table}[t]
\caption{All nine settings: mean queries for S-Track and P-Track;
$q$ is the ratio of mean P-Track to P-One queries, and $t$ is the ratio
of mean P-One to P-Track steps. C and M denote clustered and mixed blocks.
$V$ measures variation among component gradient changes, as defined in the
text; its entries are averaged over seeds.}
\label{tab:controlled-partition}
\centering\small
\begin{tabular}{rrcrrrrr}
\toprule
$n$ & $k$ & Blocks & S-Track & P-Track & $q$ & $t$ & $10^3V$\\
\midrule
128 & 64 & C & 731.5 & 1{,}439.5 & 1.039 & 4.8 & 5.778\\
512 & 64 & C & 1{,}115.5 & 4{,}861.2 & 1.050 & 15.0 & 5.778\\
2048 & 64 & C & 2{,}651.5 & 19{,}866.5 & 0.783 & 122.5 & 5.778\\
8192 & 64 & C & 8{,}795.5 & 78{,}083.2 & 0.866 & 399.4 & 5.778\\
2048 & 16 & M & 2{,}203.0 & 22{,}822.0 & 0.778 & 228.7 & 0.037\\
2048 & 16 & C & 2{,}220.2 & 24{,}956.0 & 0.857 & 213.3 & 5.657\\
2048 & 64 & M & 2{,}587.5 & 19{,}077.8 & 0.804 & 115.8 & 0.159\\
2048 & 256 & M & 4{,}126.5 & 16{,}131.8 & 0.839 & 49.9 & 0.662\\
2048 & 256 & C & 3{,}985.0 & 17{,}785.0 & 0.829 & 54.8 & 6.282\\
\bottomrule
\end{tabular}
\end{table}

\Needspace{8\baselineskip}
For the clustered setting with $n=2048$ and $k=64$, batching allows the
stepsize to increase from $1/256$ to $1/3$ and reduces mean steps from
6551.1 to 53.5. It reduces mean steps in every setting, by factors of
4.8--399.4, and mean queries in seven settings. Queries instead increase
by 3.9\% and 5.0\% at $n=128$ and $n=512$.
The theorem permits larger steps while retaining the same order of query
bound. The experiments show that the number of queries can nevertheless
increase on some instances.

In these analytic, strongly convex instances, the initial scan identifies
all exceptional components. The results illustrate optimization after this
identification and the effect of dividing an objective among components.
They do not test the worst-case cost of finding exceptions in the broader
smooth class, noisy gradients, data from applications, or the running time
of a distributed implementation.
\section{Conclusion}
We characterize the query cost of finite-sum optimization when at most
$k$ residual gradients are nonzero functions and their total oscillation
is bounded. Unknown exceptional indices contribute the additional term
$\min\{nk,\mathcal P\}$ to the optimization cost $\mathcal R$.
A common potential controls changes of the partial objective, and an
adaptive working-set method chooses its batch size from the components
already discovered. Its bound depends on the actual support size and
number of discovery stages, without requiring either quantity as input.
The lower bounds permit arbitrary adaptive queries. Interrupted searches
handle parallel chains of any positive length, while bounded sinusoidal
examples extend the matching-order result to
$\epsilon^2\leq\delta F_0\sqrt{k}/(4n)$.

The delayed-discovery experiments give every competitive method access
to the support information that its queries acquire. Adaptive batches reduce queries relative to fixed batches in every setting of the
final comparison. Their advantage over exact gradients appears when the residual
scale is sufficiently small relative to partial-objective smoothness;
full-step gradient descent wins several other settings. Population sampling
costs less in the dense setting and uses less private residual work. The separate regression study
illustrates optimization after the indices are identified. These resource
comparisons do not imply faster distributed execution.

The complexity characterization assumes exact residual gradients, free
predictor evaluations, and $\Omega\leq2F_0$. The lower-bound dimension
may grow with the parameters, and its explicit constants do not determine
numerical crossover points. The delayed examples are one-dimensional and
have structured partial objectives; their sharper smoothness bound is
supplied. Questions remaining include approximate or noisy residuals,
unknown partial smoothness, dependence on $\Omega$ independently of $F_0$,
and the cost of evaluating the predictors.

\section*{Declarations}
Generative AI tools assisted literature search, formulation and revision of
mathematical arguments, implementation of research checks, and manuscript
editing. The authors assume responsibility for all content.

\Needspace{8\baselineskip}

\end{document}